\documentclass[12pt,showkeys]{amsart}
\usepackage{color}
\usepackage{graphicx}
\usepackage{subfigure}
\usepackage{amssymb}
\usepackage{amsmath}
\usepackage{multirow}
\usepackage{enumerate}
\usepackage{epstopdf}
\usepackage{cite,stmaryrd,txfonts}

\usepackage{empheq}

\usepackage{cases}

\def\cequiv{\raisebox{-1.5mm}{$\;\stackrel{\raisebox{-3.9mm}{=}}{{\sim}}\;$}}

\newtheorem{theorem}{Theorem}
\newtheorem{remark}[theorem]{Remark}
\newtheorem{proposition}[theorem]{Proposition}
\newtheorem{lemma}[theorem]{Lemma}

\newtheorem{definition}[theorem]{Definition}

\newcounter{mnote}
\let\oldmarginpar\marginpar
\renewcommand\marginpar[1]{\-\oldmarginpar[\raggedleft\footnotesize #1]%
  {\raggedright\footnotesize #1}}

\begin{document}

\title[Optimal solvers for fourth-order PDEs]{Optimal solvers for
  fourth-order PDEs \\discretized on unstructured grids}

\author{Shuo Zhang}
\address{LSEC, Institute of Computational Mathematics and Scientific/Engineering Computing, Academy of Mathematics and System Sciences, Chinese Academy of Sciences, Beijing 100190, People's Republic of China}
\email{szhang@lsec.cc.ac.cn}
\author{Jinchao Xu}
\address{Center for Computational Mathematics and Applications, Department of Mathematics, The Pennsylvania State University, University Park, Pennsylvania 16802}
\email{xu@math.psu.edu}

\subjclass[2000]{Primary 65F10, 65N22, 65N30, 46E35}
\keywords{Fourth-order problem, unstructured grid, finite element, optimal solver, preconditioning, Sobolev space}

\begin{abstract}
This paper provides the first provable $\mathcal{O}(N \log N)$ algorithms for the linear system arising from the direct finite element discretization of the fourth-order equation with different boundary conditions on unstructured grids of size $N$ on an arbitrary polygoanl domain. Several preconditioners are presented, and the conjugate gradient methods applied with these preconditioners are proven to converge uniformly with respect to the size of the preconditioned linear system. One main ingredient of the optimal preconditioners is a mixed-form discretization of the fourth-order problem. Such a mixed-form discretization leads to a non-desirable ---either non-optimal or non-convergent--- approximation of the original solution, but it provides optimal preconditioners for the direct finite element problem. It is further shown that the implementation of the preconditioners can be reduced to the solution of several discrete Poisson equations. Therefore, any existing optimal or nearly optimal solver, such as geometric or algebraic multigrid methods, for Poisson equations would lead to a nearly optimal solver for the discrete fourth-order system. A number of nonstandard Sobolev spaces and their discretizations defined on the boundary of polygonal domains are carefully studied and used for the analysis of those preconditioners.

\end{abstract}

\maketitle

\tableofcontents

\section{Introduction}
In this paper, we study numerical methods for solving finite element systems for boundary value problems for fourth-order elliptic partial differential equations on unstructured grids of an arbitrary polygonal domain. The fourth-order partial differential equation has applications in, for example, solid and fluid mechanics and material sciences. Many different finite element methods, including conforming and nonconforming, have been developed in the literature for discretizing the boundary value problems.  All these discretizations lead to very ill-conditioned linear systems with condition numbers of order $\mathcal{O}(h^{-4})$ that are difficult to solve.

Multigrid methods are among the most efficient techniques for solving these systems.  In particular, geometric multigrid methods based on a nested sequence of multilevel geometric grids have been extensively studied in the literature for fourth-order problems, c.f. \cite{BrambleZhang1995, Brenner1989, Brenner1999, Hanisch1993, PeiskerBraess1987, Stevenson2003, Xu2002, Zhang1989, Zhao2004} and references therein. The efficiency of these methods, however, depends crucially on appropriate underlying multilevel structures.  Because such multilevel structures are not naturally available in most unstructured grids in practice, multigrid methods of this type are generally quite difficult to use.  More user-friendly methods such as algebraic multigrid methods (that can be applied to unstructured grids) have also been studied in the literature, see \cite{BrandtMcCormickRuge1982, Brandt1986, RugeStuben1987, VanekMandelBrezina1996, Livne2004}.  But the efficiency of these methods applied on the fourth-order finite element problem is limited. Further, there is still no theory to support the methods of these types. In any event, to the authors' knowledge, no mathematically optimal solvers for fourth-order finite element problems discretized on unstructured grids are presented in the literature.

In this paper, we develop a class of methods that fall in between geometric and algebraic multigrid methods. We use the final geometric grid only which has no hierarchy and can be entirely unstructured. We circumvent the obstacles inhering in a lack of natural hierachical structure by exploring the deep virtue of multilevel methodology in the framework of the Fast Auxiliary Space Preconditiong (FASP) method that exists in the literature \cite{Xu1996, Xu2010}. These methods combine the practical advantage of algebraic multigrid methods in that they are easy to use with the solid theoretical foundation associated with geometric multigrid methods.

To cope with the conforming and nonconforming finite element
discretizations for the fourth-order problem with different boundary conditions, we present several preconditioners and then use the conjugate gradient method on the preconditioned linear systems. One main ingredient in the construction of these preconditioners is an auxiliary discretization of a mixed form of the fourth-order problem whereby piecewise linear finite elements are used for both the original variable and the auxiliary variable. Such  simple mixed finite element discretizations often lead to non-desirable (either non-optimal or non-convergent) approximations of the original solution; however, a proper combination of the solution and some elementary point relaxation iterative methods for the original system, such as the Jacobi and symmetric Gauss--Seidel methods, produces a preconditioner that can capture the spectrum of the original linear systems well.  And, as a result, the conjugate gradient method applied to the preconditoned system converges uniformly with respect to the size of the systems. The solution of the linear mixed system will be reduced to the solution of several discrete Poisson equations. Thus, any existing optimal or nearly optimal Poisson solver, such as geometric or algebraic multigrid methods, would lead to a nearly optimal solver for the discrete fourth-order system. These solvers work for a relatively large class of finite element systems. 

In order to analyze the optimality and complexity of the preconditioners, a number of technical results associated with the Sobolev space $H^2\cap H_0^1(\Omega)$ on the polygonal domain are developed in the paper. That is, the finite element discretizations and the continuous and discrete traces on the boundary are studied, as the generalization of some existing results for the smooth domain. Based on the properties of the newly developed trace spaces and the trace and extension operators defined on the trace spaces, the optimality of each solver is shown rigorously on both convex and nonconvex domains. 

With various numerical examples also reported to support the theoretical results, this paper appears to be the first to present provable $\mathcal{O}(N\log N)$ algorithms for the fourth-order equation discretized on unstructured grids.

The rest of the paper is organized as follows. In the next section,
we introduce the model problem and basic definitions and
notations of finite element methods. In Section \ref{sec:pre}, we introduce some preliminaries of the paper, namely a concise introduction to preconditioned conjugate gradient (PCG) method and FASP. In Section \ref{sec:dislaps}, we discuss several discrete second-order operators defined on finite element spaces, and in Section \ref{sec:reg}, we develop some trace spaces of specific Sobolev spaces and discuss their discretization. Sections \ref{sec:opsol} and \ref{sec:numexp} are devoted to the optimal solver for the finite element problem of the fourth-order problems. Both theoretical results and numerical examples are provided to verify the optimality of the solvers. Finally in Section \ref{sec:con}, conclusions and remarks are given.

\section{Model problems and finite element discretizations}\label{sec:model}
\subsection{Model problems}
Let $\Omega\subset\mathbb{R}^2$ be a bounded polygon, with $\Gamma=\partial\Omega$ its boundary, and $\Gamma_i$ and $i=1,\dots,K$ as the edges of $\Omega$, such that $\Gamma=\bigcup_{i=1}^K \Gamma_i$. We consider the model biharmonic equation
\begin{equation}
\Delta^2u=f
\end{equation}
on $\Omega$ equipped with boundary value conditions of the first and second kind, and its variational formulation. In doing so, we find $u\in M_k$, such that
\begin{equation}\label{eqn:modelsd}
(\nabla^2u:\nabla^2v)=(f,v)\quad\forall\,v\in M_k,
\end{equation}
where $M_k$ is a Hilbert space with a certain boundary value condition; namely, $M_1=H^2_0(\Omega)$ when the Dirichlet boundary condition of the first kind is considered, and $M_2=H^2(\Omega)\cap H^1_0(\Omega)$, when the Dirichlet boundary condition of the second kind is considered. In terms of elasticity, the first is in accordance with the clamped plate and the second is in accordance with the simply supported plate, and they will be referred to thereafter as the first and second biharmonic problem, respectively. It is known for $w\in H^2(\Omega)\cap H^1_0(\Omega)$ that $\|\Delta w\|_{0,\Omega}=\|\nabla^2w\|_{0,\Omega}=|w|_{2,\Omega}\cequiv \|w\|_{2,\Omega}$(\cite{Grisvard1992}). The well-posedness of the boundary value problem is obvious by the Lax--Milgram Lemma. An equivalent weak form of the boundary value problem is to seek $u\in M_k$, such that
\begin{equation}
(\Delta u,\Delta v)=(f,v)\quad\forall\,v\in M_k.
\end{equation}

\subsection{Finite element problems for model problems}

\subsubsection{Triangulation and finite element}\label{sec:pre:fes}

Let $\mathcal{T}_h$ be a quasi-uniform triangular triangulation of domain $\Omega$, $\Omega=\cup_{T\in\mathcal{T}_h}T$. We denote $h_T$ as the mesh size of $T$,  and $h$ as the meshsize. 

For a given triangulation, the finite element space is defined by the elementwise shape function space and the continuity of the nodal parameters. In this paper, we study the finite elements whose nodal parameters are of the type
$\mathcal{N}_\alpha(\varphi)=\fint_{D_\alpha} (\nabla^{k_\alpha}\varphi)(\mathbf{t}_1,\dots,\mathbf{t}_{k_\alpha})$,
where $D_\alpha$ is the integral domain with respect to $\alpha$, and $\mathbf{t}_1,\dots,\mathbf{t}_{k_\alpha}$ are $k_\alpha$ identical or different unit vectors to denote the direction of the derivative, $k_\alpha$=0,\,1,\,2,\dots. In general, $D_\alpha$ is a subsimplex of the triangulation; that is, it can be a point, an edge, or a triangle. When $D_\alpha$ is a point, the average of integration is reduced to the evaluation on the point. Define the degree of the nodal parameter by $deg(\mathcal{N}_\alpha)=deg(\alpha):=k_\alpha.$ When $deg(\alpha)=1$, there is only one direction involved for the derivative, and the direction is denoted by $\mathbf{t}_\alpha$.

Let $K$ be a subsimplex of the triangulation. Define the neighboring patch of $K$ as $\omega_K=\bigcup\{T'\in\mathcal{T}_h: \overline{T'}\cap\overline{K}\neq\emptyset\}.$ Given $T\in\mathcal{T}_h$, let $P_k(T)$ denote the polynomials on $T$ with a degree not higher than $k$, and $P_k(\mathcal{T}_h)$ the set of piecewise polynomials that belong to $P_k(T)$ at each element $T$.

For a given finite element space, the nodal basis is defined to be the set of all the dual basis functions with respect to the nodal parameters. Let $\varphi_\alpha$ be the nodal basis function with respect to $\mathcal{N}_\alpha$. Define $\omega_\alpha:=\bigcup\{T:\mathring{T}\bigcap supp(\varphi_\alpha)\neq\emptyset\}$, $\#\omega_\alpha:=\#\{T:\mathring{T}\bigcap supp(\varphi_\alpha)\neq\emptyset\}$, and $h_\alpha=\max_{T\subset\omega_\alpha}h_T$. Any $\omega_\alpha$ can be written as $\omega_\alpha=T_1\bigcup\dots\bigcup T_{\#\omega_\alpha},$ such that $T_j$ and $T_{j+1}$ share a common edge.

We make use of the following notations. Let $\mathcal{N}_h$ denote the set of all the vertices, $\mathcal{N}_h=\mathcal{N}_h^i\cup\mathcal{N}_h^b$, with $\mathcal{N}_h^i$ and $\mathcal{N}_h^b$ consisting of the interior vertices and the boundary vertices, respectively. Similarly, let $\mathcal{E}_h=\mathcal{E}_h^i\bigcup\mathcal{E}_h^b$ denote the set of all the edges, with $\mathcal{E}_h^i$ and $\mathcal{E}_h^b$ consisting of the interior edges and boundary edges, respectively. For an edge $e$, $\mathbf{n}_e$ is a unit vector normal to $e$. For $e\in\mathcal{E}_h^i$, $T_e^L$ and $T_e^R$ are the two adjacent elements that share $e$ as the common edge, and $\mathbf{n}_e^L$ and $\mathbf{n}_e^R$ denote the unit outer normal vectors of $T_e^L$ and $T_e^R$ on $e$. Define on the edge $e$
$$
\left\llbracket\frac{\partial \psi}{\partial\mathbf{n}_e}\right\rrbracket:=\left\{\begin{array}{l}\displaystyle\frac{\partial\psi|_{T_e^L}}{\partial\mathbf{n}_e^L}+\frac{\partial\psi|_{T_e^R}}{\partial\mathbf{n}_e^R}=\frac{\partial\psi|_{T_e^L}}{\partial\mathbf{n}_e}(\mathbf{n}_e\cdot\mathbf{n}_e^L)+\frac{\partial\psi|_{T_e^R}}{\partial\mathbf{n}_e}(\mathbf{n}_e\cdot\mathbf{n}_e^R),\\ \displaystyle \qquad \qquad\mbox{when}\ e\in\mathcal{E}_h^i\ \mbox{and}\ \ e=T_e^L\cap T_e^R;\\
\displaystyle\frac{\partial \psi}{\partial\mathbf{n}_e},\quad\quad\mbox{when}\ e\in\mathcal{E}_h^b.
\end{array}\right.
$$

\subsubsection{Finite element problems}\label{sec:fempro}

We assume that the finite element space under consideration satisfies the following conditions:
\begin{description}
\item[A1] the quadratic polynomials are contained in the elementwise shape function space;
\item[A2] each $w_h\in M_h$ is weakly continuous and $\nabla w_h$ is weakly continuous on each $e\in\mathcal{E}_h^i$ in the following sense, where $\mathbf{n}$ is the normal vector of $e$ and $\tau$ is the tangential normal vector of $e$:  $w_h$ is continuous on at least one point on $e$, $\partial_\mathbf{n} w_h$ is continuous on at least one point on $e$, and $\partial_\tau w_h$ is continuous on at least one point on $e$;
\item[A3] for $w_h\in M_h$ and $T\in\mathcal{T}_h$,
\begin{equation}
\|w_h\|_{0,T}^2\cequiv\sum_{\omega_\alpha\supset T}h_T^{2+2deg(\alpha)}\mathcal{N}_\alpha(w_h)^2.
\end{equation}
\end{description}
Here and after, we make use of $\lesssim$, $\gtrsim$ and $\cequiv$ to denote $\leqslant$, $\geqslant$ and $=$ up to a constant. The hidden constants depend on the domain. And, when the triangulation is involved, they also depend on the shape-regularity of the triangulation, but they do not depend on $h$ or any other mesh parameter.
\begin{remark}
The Assumptions A1--A3 are mild, and they hold for most finite elements for fourth-order problems, including the Morley element, the modified-Zienkiewicz element\cite{WangShiXu2007}, the Nilssen--Tai--Winther element\cite{NilssenTaiWinther2001}, the Morley--Zienkiewicz element\cite{ShiWangContext}, the Bell element, and the Argyris element.
\end{remark}

We define a piecewise $H^2$ function space as
$$
H^2(\mathcal{T}_h):=\{w\in L^2(\Omega):w|_T\in H^2(T),\ \forall\,T\in\mathcal{T}_h\}.
$$
The operator $\nabla^2_h$ is defined for $w\in H^2(\mathcal{T}_h)$ by
$(\nabla^2_h w)|_T=\nabla^2(w|_T).$ By the assumptions A1 and A2, the functional $|\cdot|_{2,h}$ defined as $|w_h|_{2,h}=\|\nabla_h^2w_h\|_{0,\Omega}$ is a nontrivial semi-norm on $M_h$.

In practice, the finite element space $M_{h,1}$ that is associated with $H^2_0(\Omega)$ is defined by
$$
\quad M_{h,1}=\Big\{w_h\in M_h:\ \mathcal{N}_\alpha(w_h)=0, \mbox{when}\ deg(\alpha)\leqslant\ 1\ \mbox{and}\ D_\alpha\subset\partial\Omega\Big\}.\qquad
$$ 
And the finite element space $M_{h,2}$ that is associated with $H^2(\Omega)\cap H^1_0(\Omega)$ is defined by
$$
M_{h,2}:=\Bigg\{w_h\in M_h:\ \displaystyle \fint_{D_\alpha}(\partial_\tau)^{deg(\alpha)} w_h=0\ \mbox{when}\ D_\alpha\subset\partial\Omega,\ \mbox{and}\  deg(\alpha)=0,1\Bigg\}.
$$
By definition, $M_{h,1}\subset M_{h,2}$. The semi-norm $|\cdot|_{2,h}$ is a norm on $M_{h,2}$.

The finite element discretization of \eqref{eqn:modelsd} is to find $u_h\in M_{h,k}$ ($k=1,2$) such that
\begin{equation}\label{vpnc}
a_h(u_h,v_h):=(\nabla^2_h u_h:\nabla^2_hv_h)=\int_\Omega fv_h,\quad\forall\,v_h\in M_{h,k}.
\end{equation}
We introduce the operator $A_{h,k}:M_{h,k}\to M_{h,k}$,
$$
(A_{h,k} v_h,w_h)=a_h(v_h,w_h),\quad\forall\,v_h,w_h\in M_{h,k}.
$$
Then the linear system to solve is $A_{h,k} u_h=f_h$, with $f_h$ the $L^2$ projection of $f$ on $M_{h,k}$.

\section{Preliminaries: PCG and FASP}\label{sec:pre}
\subsection{Preconditioned conjugate gradient (PCG) method}\label{sec:pre:pcg}

The PCG method is the basis of all the preconditioning techniques studied in this paper. For a linear system $Au=f$ in which $A$ is a symmetric positive definite (SPD) operator, the PCG method can be viewed as a conjugate gradient method applied to the preconditioned system $BAu = Bf$. Here, $B$ is an SPD operator, and $BA$ is symmetric with respect to the inner product $(\cdot,\cdot)_{B^{-1}}:=(B^{-1}\cdot,\cdot)$.

Let $u_k$, $k=0,1,2,\cdots$ be the solution sequence of the PCG algorithm. It is well known that
\begin{equation}\label{xz2.3}
\|u-u_k\|_A \leqslant 2 \left(\frac{\sqrt{\kappa(BA)}-1}{\sqrt{\kappa(BA)}+1}\right)^k \|u-u_0\|_A,
\end{equation}
which implies that generally the smaller the condition number $\kappa(BA)$, the faster the PCG method will converge. However, though it is sufficient for many applications, the estimate given in \eqref{xz2.3} is not sharp. An improved estimate can be obtained in terms of the eigenvalue distribution of $BA$ \cite{Axelsson2003, Hackbusch1994, XuPSU}. More specifically, we focus on the case in which we can divide $\sigma(BA)$, the spectrum of $BA$, into two subsets, $\sigma_0(BA)$ and $\sigma_1(BA)$, where the eigenvalues in $\sigma_1$ are bounded from above and from below and $\sigma_0$ consists of all the ``bad" eigenvalues remaining. The lemma below can be found in many references, e.g. \cite{Axelsson2003}.

\begin{lemma}\label{lem:pcgrate}
Suppose that $\sigma(BA) = \sigma_0 (BA)\cup \sigma_1 (BA)$ such that there are $m$ elements in $\sigma_0 (BA)$ and $\lambda\in [a, b]$ for each $\lambda\in \sigma_1 (BA)$. Then
\begin{equation}\label{xz2.4}
\|u-u_k\|_A\leqslant 2K\Bigg(\frac{\sqrt{b/a}-1}{\sqrt{b/a}+1}\Bigg)^{k-m}\|u-u_0\|_A
\end{equation}
where
$$
K = \max_{\lambda\in \sigma_1 (BA)} \prod_{\mu\in\sigma_0(BA)}|1-\frac{\lambda}{\mu}|.
$$
\end{lemma}

If the eigenvalues of $BA$ are ordered as $0 < \lambda_1(BA)\leqslant \cdots\leqslant \lambda_{n-m_b}(BA)\ll \lambda_{n-m_b+1}(BA)\leqslant \cdots\leqslant \lambda_n(BA)$, then the convergence rate estimate \eqref{xz2.4} becomes
\begin{equation}\label{xz2.5}
\frac{\|u-u_k\|_A}{\|u-u_0\|_A}\leqslant 2\left(\frac{\sqrt{\lambda_{n-m_b}(BA)/\lambda_1(BA)}-1}{\sqrt{\lambda_{n-m_b}(BA)/\lambda_1(BA)}+1}\right)^{k-m_b}, \quad k\geqslant m_b.
\end{equation}
In this case, $\lambda_{n-m_b}(BA)/\lambda_1(BA)$ plays a dominant role in determining the convergence rate, and it is sometimes called an effective condition number of BA. In general, let $A$ be an SPD operator, and the eigenvalues of $A$ be ordered by $0<\lambda_1(A)\leqslant \lambda_2(A)\leqslant\dots\leqslant\lambda_n(A)$. Define the $m-$th effective condition number of $A$ by
$\kappa^{\rm eff}_m(A):=\lambda_{n-m}(A)/\lambda_1(A).$
A similar consideration can be found in \cite{XuZhu2008}. Obviously, $\kappa^{\rm eff}_0(A)=\kappa(A)$ and $\kappa^{\rm eff}_m(A)$ decreases as $m$ grows. Then based on \eqref{xz2.5}, given a tolerance $\varepsilon < 1$, the number of iterations of the PCG algorithm needed for $\frac{\|u-u_k\|_A}{\|u-u_0\|_A}\leqslant\varepsilon$ is
\begin{equation}\label{eq:effcondest}
k\gtrsim m_b+\log\varepsilon(\kappa^{\rm eff}_{m_b}(A))^{1/2}.
\end{equation}

To estimate the effective condition numbers, a basic tool is the Courant minimax principle (see e.g., \cite{GolubLoan1996}), which leads to the following lemma:
\begin{lemma}\label{lem:minmax}
Let $V$ be a Hilbert space and $A : V \to V$ be an SPD operator on $V$. Suppose the eigenvalues of $A$ are ordered as $ \lambda_1\leqslant\lambda_2\leqslant\cdots\leqslant\lambda_n$, then for any subspace $V_0\subset V$ with $dim(V_0)=m$, the following estimates hold:
\begin{equation}\label{xz2.7}
\lambda_m\leqslant\max_{v\in V_0}\frac{(Av,v)}{(v,v)}\quad\mbox{and}\quad
\lambda_{n+1-m}(A)\geqslant\min_{v\in V_0}\frac{(Av,v)}{(v,v)}.
\end{equation}
\end{lemma}
\subsection{Theory of the Fast Auxiliary Space Preconditioning (FASP) method}
\label{sec:pre:fasp}

In this section, we give a summary of preconditioning techniques based on auxiliary spaces as developed in \cite{Xu1992, Xu1996, HiptmairXu2007, Xu2010}.

Let $V$ stand for a real Hilbert space with an inner product $a(\cdot,\cdot)$ and a (energy) norm $\|\cdot\|_A$. Write $A : V \to V'$ for the isomorphisms associated with $a(\cdot,\cdot)$. Let $s(\cdot,\cdot)$ be another inner product on $V$, and $S: V\mapsto V'$ the isomorphism associated with $s(\cdot,\cdot)$. Let $W_1,\dots,W_J$, $J\in\mathbb{N}$ be Hilbert spaces endowed with inner products $\bar{a}_j(\cdot,\cdot)$, $j=1,\dots,J$. Write $A_j: W_j\mapsto W_j'$ for the isomorphisms associated with $a_j(\cdot,\cdot)$, $j=1,\dots,J$. Furthermore, for each $W_j$, we need a linear transfer operator $\Pi_j: W_j \mapsto V$. We tag the adjoint operators by $``{}^*"$.

Then the fast auxiliary space preconditioner (see \cite{Xu1996,Xu2010}) is defined by
\begin{equation}\label{asp}
B=R+\sum_{j=1}^J\Pi_j\circ A_j^{-1}\circ\Pi_j^*,\ \ \mbox{with}\ R=S^{-1}.
\end{equation}

\begin{theorem}\label{thm:fasp}
\cite{Xu1996, HiptmairXu2007} Assume that
\begin{enumerate}
\item There are constants $c_j > 0$, such that
\begin{equation*}\label{assum1}
\|\Pi_jw_j\|_A\leqslant c_ja_j(w_j,w_j)^{1/2},\, \forall\, w_j\in W_j;
\end{equation*}
\item There is a constant $c_s>0$, such that
$$
\|v\|_A\leqslant c_ss(v,v)^{1/2},\  \ \forall\,v\in V;
$$
\item There is a constant $c_0>0$, such that for any $v\in V$, there are $v_0\in V$ and $w_j\in W_j$ such that
$$
v=v_0+\sum_{j=1}^J\Pi_jw_j,
$$
and
\begin{equation*}
s(v_0,v_0)+\sum_{j=1}^Ja_j(w_j,w_j)\leqslant c_0^2\|v\|_A^2.
\end{equation*}
\end{enumerate}
Then
$$\kappa(BA)\leqslant c_0^2\Big(c_s^2+c_1^2+\cdots+c_J^2\Big).$$
\end{theorem}

In applications, $V$ and all $W_j$ are usually finite element spaces with bases consisting of locally supported functions. Plugging basis functions
into the bilinear forms would lead to the algebraic representation of the preconditioner.

\section{Discrete Second-order operators}\label{sec:dislaps}

The Laplacian operator and the Hessian operator are second-order operators of fundamental importance. In this section, we study their discretization on finite element spaces. 

For a linear space $M$ and its subspace $M_*\subset M$, denote $codim(M_*,M)$ as the codimension of $M_*$ in $M$. By definition, we can obtain the following lemma directly.
\begin{lemma}\label{lem:codim}
Assume a linear space $M$ and a linear operator $L$ on $M$. For any subspace $M_*$ of $M$, $codim(L(M_*),L(M))\leqslant codim(M_*,M)$; if $ker(L)\subset M_*$, then $codim(L(M_*),L(M)) = codim(M_*,M)$.
\end{lemma}

Assume $\mathcal{R}(\Delta)=\Delta(H^2(\Omega)\cap H^1_0(\Omega)):=\{\Delta w:w\in H^2(\Omega)\cap H^1_0(\Omega)\}$.  Following, e.g., \cite{Grisvard1992}, we have $codim(\mathcal{R}(\Delta),L^2(\Omega))=m_0$, where $m_0$ is the number of reentrant corners on $\Gamma$ throughout this paper, and by a ``reentrant corner" we refer to the corner whose interior angle is bigger than $\pi$. When the domain $\Omega$ is convex, $m_0=0$.

\subsection{Discrete Laplacians on linear finite element spaces for $H^1(\Omega)$}\label{sec:les:dle}

Given a triangulation $\mathcal{T}_h$, denote by $V_h\subset H^1(\Omega)$ the continuous linear element space on $\mathcal{T}_h$, and $V_{h0}=V_h\cap H^1_0(\Omega)$. Moreover, define $V_{hb}$ as the complementary subspace of $V_{h0}$ in $V_h$, namely, $V_{hb}\subset V_h\ \mbox{and}\ V_h=V_{h0}\oplus V_{hb}$, and $B_h$ as the trace space of $V_h$ on $\Gamma$. Define $\mathring{\mathcal{B}}_h$ as consisting of functions that vanish on the corners of $\Gamma$, and $\mathcal{B}_h^c$ the complementary of $\mathring{\mathcal{B}}_h$ in $\mathcal{B}_h$.

\begin{definition}\label{def:2ndoprs}
Define the discrete Laplacian operators of first and second kind on $V_{h0}$ as
\begin{equation}
-\Delta_{h,k}:\ V_{h0}\to V_{h,k}\ \mbox{such\ that}\ (-\Delta_{h,k} w_h,v_h)=(\nabla w_h,\nabla v_h),\quad\forall\,v_h\in V_{h,k},\ \forall\,w_h\in V_{h0};
\end{equation}
where $V_{h,1}=V_h$, and $V_{h,2}=V_{h0}$.
\end{definition}

Define $\widetilde{V}_{h0}:=\Big\{p_h\in V_{h0}:\Delta_{h,2}p_h\in\mathcal{R}(\Delta)\Big\}$, then $\Delta_{h,2}\widetilde{V}_{h0}=\mathcal{R}(\Delta)\cap V_{h0}$. As $ker(\Delta_{h,2})\subset \widetilde{V}_{h0}$, $codim(\widetilde{V}_{h0},V_{h0})=codim(\Delta_{h,2}\widetilde{V}_{h0},\Delta_{h,2}V_{h0})= codim(\mathcal{R}(\Delta)\cap V_{h0},V_{h0})\leqslant m_0$.

Each of the two Laplacian operators has a second-order equivalent description in terms of normal derivative jumps on element edges. Let $\mathcal{E}_{h,1}=\mathcal{E}_{h}$, and $\mathcal{E}_{h,2}=\mathcal{E}^i_{h}$.

\begin{lemma}\label{lem:dlevc}
The Laplacian operators can be described by
\begin{equation}\label{controltwoside}
\|\Delta_{h,k}p_h\|_{0,\Omega}^2
\lesssim
\sum_{e\in\mathcal{E}_{h,k}}h_e^{-1}\int_e\Big\llbracket\frac{\partial p_h}{\partial\mathbf{n}_e}\Big\rrbracket^2\lesssim
\|\Delta_{h,k}p_h\|_{0,\Omega}^2,
\end{equation}
which holds for all $p_h\in V_{h0}$ and $k=1,2$, with this exception: when $k=2$,  the right inequality above only holds for $p_h\in \widetilde{V}_{h0}$.
\end{lemma}
\begin{proof} Let $p_h\in V_{h0}$ and $k=1,2$, then
$$
\begin{array}{rl}
\displaystyle\|\Delta_{h,k}p_h\|_{0,\Omega}^2&\displaystyle=-(\nabla p_h,\nabla \Delta_{h,k}p_h) =-\sum_{e\in\mathcal{E}_{h,k}}\int_e\Big\llbracket\frac{\partial p_h}{\partial\mathbf{n}_e}\Big\rrbracket \Delta_{h,k}p_h\\
&\displaystyle=-\sum_{e\in\mathcal{E}_{h,k}}\Big(h_e^{-1}\int_e\Big\llbracket\frac{\partial p_h}{\partial\mathbf{n}_e}\Big\rrbracket\Big)\Big(\int_e\Delta_{h,k}p_h\Big)\\
&\displaystyle\leqslant \Big(\sum_{e\in\mathcal{E}_{h,k}}h_e^{-2}\Big(\int_e \Big\llbracket\frac{\partial p_h}{\partial\mathbf{n}_e}\Big\rrbracket\Big)^2\Big)^{1/2}\Big(\sum_{e\in\mathcal{E}_{h,k}}\big(\int_e\Delta_{h,k}p_h\big)^2\Big)^{1/2}\\
&\displaystyle\cequiv (\sum_{e\in\mathcal{E}_{h,k}}h_e^{-1}\int_e \Big\llbracket\frac{\partial p_h}{\partial\mathbf{n}_e}\Big\rrbracket^2)^{1/2}\|\Delta_{h,k}p_h\|_{0,\Omega}.
\end{array}
$$
Namely, $\displaystyle\|\Delta_{h,k}p_h\|_{0,\Omega}^2 \lesssim \sum_{e\in\mathcal{E}_{h,k}}h_e^{-1}\int_e\Big\llbracket\frac{\partial p_h}{\partial\mathbf{n}_e}\Big\rrbracket^2.$

Now we turn to the right inequality of \eqref{controltwoside}. We first consider the case $k=2$. Let $p\in H^1_0(\Omega)$ be the unique solution of
$$
\left\{
\begin{array}{ll}
\displaystyle-\Delta p=-\Delta_{h,2}p_h,&\mbox{in}\,\Omega,\\
\displaystyle p=0,&\mbox{on}\,\partial\Omega.
\end{array}
\right.
$$
Then, when $\Delta_{h,2}p_h\in\mathcal{R}(\Delta)$, it holds that $|p_h-p|_{1,\Omega}\lesssim h\|p\|_{2,\Omega}\lesssim h\|\Delta_{h,2}p_h\|_{0,\Omega}.$
Thus,
\begin{equation}\label{convex}
\begin{array}{rl}
\displaystyle\sum_{e\in\mathcal{E}_{h,2}}h_e^{-1} \int_e\Big\llbracket\frac{\partial p_h}{\partial\mathbf{n}_e}\Big\rrbracket^2&\displaystyle=\sum_{e\in\mathcal{E}_{h,2}}h_e^{-1}\int_e \Big\llbracket\frac{\partial (p_h-p)}{\partial\mathbf{n}_e}\Big\rrbracket^2\\
&\displaystyle\lesssim\sum_{e\in\mathcal{E}_{h,2}}\sum_{\partial T\supset F}\Big(h_e^{-2}|p-p_h|_{1,T}^2+|p-p_h|_{2,T}^2\Big)\\
&\displaystyle\lesssim |p|_{2,\Omega}^2\lesssim \|\Delta_{h,2}p_h\|_{0,\Omega}^2.
\end{array}
\end{equation}
The right inequality of \eqref{controltwoside} is proved for $k=2$, and we turn to the case $k=1$. As
$$
\|\Delta_{h,1}p_h\|_{0,\Omega}=\sup_{q_h\in V_h\setminus\{0\}}\frac{(\Delta_{h,1} p_h,q_h)}{\|q_h\|_{0,\Omega}}=\sup_{q_h\in V_h\setminus\{0\}}\frac{(\nabla p_h,\nabla q_h)}{\|q_h\|_{0,\Omega}},
$$
we will find some $r_h\in V_h\setminus\{0\}$, such that
\begin{equation}\label{supre}
\sum_{e\in\mathcal{E}_{h,1}}h_e^{-1}\int_e \Big\llbracket\frac{\partial p_h}{\partial\mathbf{n}_e}\Big\rrbracket^2 \lesssim\frac{(\nabla p_h,\nabla r_h)^2}{\|r_h\|_{0,\Omega}^2}.
\end{equation}

If $\Omega$ is convex, then by the same argument as for $k=2$ with the homogeneous Dirichlet problem replaced by the homogeneous Neumann problem, we obtain that $\displaystyle\sum_{e\in\mathcal{E}_{h,1}}h_e^{-1}\int_e \Big\llbracket\frac{\partial p_h}{\partial\mathbf{n}_e}\Big\rrbracket^2 \lesssim\|\Delta_{h,1}p_h\|_{0,\Omega}^2$. Then \eqref{supre} holds for $r_h:=-\Delta_{h,1}p_h$.

If $\Omega$ is not convex, let $\hat{\Omega}$ be the convex hull of all the corner points of $\Omega$, and let $\hat{\mathcal{T}}_h$ be a quasi-uniform triangulation on $\hat{\Omega}$ such that $\mathcal{T}_h$ is a subtriangulation of $\hat{\mathcal{T}}_h$. Similarly, the set of edges $\hat{\mathcal{E}}_h$ is defined on $\hat{\mathcal{T}}_h$, and the discrete Laplace operator $\Delta_{\hat{h},1}$ is defined on $\hat{V}_h$, the continuous linear element space defined on $\hat{\mathcal{T}}_h$. Denote by $\hat{p}_h$ the extension of $p_h$ to $\hat{V}_h$, such that $\hat{p}_h|_{\hat{\Omega}\setminus\Omega}=0$. Then
$$
\sum_{e\in\mathcal{E}_{h,1}}h_e^{-1}\int_e \Big\llbracket\frac{\partial p_h}{\partial\mathbf{n}_e}\Big\rrbracket^2=\sum_{e\in\hat{\mathcal{E}}_h}h_e^{-1}\int_e \Big\llbracket\frac{\partial \hat{p}_h}{\partial\mathbf{n}_e}\Big\rrbracket^2.
$$
Denote $\hat{r}_h=-\Delta_{\hat{h},1}\hat{p}_h$, then $\displaystyle\sum_{e\in\hat{\mathcal{E}}_h}h_e^{-1}\int_e \Big\llbracket\frac{\partial \hat{p}_h}{\partial\mathbf{n}_e}\Big\rrbracket^2\leqslant C_1\|\hat{r}_h\|^2_{0,\hat{\Omega}}$ for a constant $C_1$ depending on $\hat{\Omega}$ and the shape regularity of $\hat{\mathcal{T}}_h$ only. Define $r_h:=\hat{r}_h|_\Omega$ so that
$$
\begin{array}{rl}
\displaystyle\|r_h\|_{0,\Omega}^2\sum_{e\in\mathcal{E}_{h,1}}h_e^{-1}\int_e \Big\llbracket\frac{\partial p_h}{\partial\mathbf{n}_e}\Big\rrbracket^2  & \displaystyle \leqslant \|\hat{r}_h\|_{0,\hat{\Omega}}^2\sum_{e\in\hat{\mathcal{E}}_h}h_e^{-1}\int_e \Big\llbracket\frac{\partial \hat{p}_h}{\partial\mathbf{n}_e}\Big\rrbracket^2\\
&\displaystyle\leqslant C_1\|\hat{r}_h\|_{0,\hat{\Omega}}^4 =\Big(\int_{\hat{\Omega}}\nabla\hat{r}_h\nabla\hat{p}_h\Big)^2=\Big(\int_\Omega\nabla r_h\nabla p_h\Big)^2.
\end{array}
$$
Therefore, \eqref{supre} holds for this $r_h$. The right inequality of \eqref{controltwoside} is proved for $k=1$, and this finishes the proof of the lemma.
\end{proof}
\subsection{Discrete Hessions on finite element spaces for $H^2(\Omega)$}\label{sec:les:fesf}
Let $M_h$ be a finite element space for $H^2(\Omega)$ as described in Section~\ref{sec:fempro}. Define interpolation operators on finite element spaces as
$$
\begin{array}{ll}
I_h:M_{h,2}\to V_{h0}, &\displaystyle (I_hw_h)(a)=\frac{1}{\#\omega_a}\Big(\sum_{T\subset\omega_a}w_h|_T(a)\Big),\ \forall\,a\in \mathcal{N}_h^i,\ \mbox{and}\ \forall\, w_h\in M_{h,2};\\
\Pi_{h,1}:\displaystyle V_{h0}\to M_{h,1}, 
&\mathcal{N}_\alpha(\Pi_{h,1}p_h)=0,\ \mbox{when}\ D_\alpha\subset\partial\Omega;\\
&\displaystyle\mathcal{N}_\alpha(\Pi_{h,1}p_h)=\frac{1}{\#\omega_\alpha}\sum_{T\subset\omega_\alpha}\mathcal{N}_\alpha(p_h|_T),\ \mbox{else};\\ 
\Pi_{h,2}:\displaystyle V_{h0}\to M_{h,2}, 
& \displaystyle \mathcal{N}_\alpha(\Pi_{h,2}p_h)= \frac{1}{\#\omega_\alpha}\sum_{T\subset\omega_\alpha}\fint_{D_\alpha}\partial_\mathbf{n}(p_h|_T)(\mathbf{n}\cdot\mathbf{t}_\alpha),\ \mbox{when}\ D_\alpha\subset\partial\Omega,\  \mbox{and}\ deg(\alpha)=1,\\
& \displaystyle\mathcal{N}_\alpha(\Pi_{h,2}p_h)=\frac{1}{\#\omega_\alpha}\sum_{T\subset\omega_\alpha}\mathcal{N}_\alpha(p_h|_T),\ \mbox{else}.
\end{array}
$$
Define $\widetilde{M}_{h,2}:=\Big\{w_h\in M_{h,2}:\Delta_{h,2}I_hw_h\in\mathcal{R}(\Delta)\Big\}$.  Then, $w_h\in\widetilde{M}_{h,2}$ if and only if $I_hw_h\in\widetilde{V}_{h0}$. By Lemma \ref{lem:codim}, $codim(\widetilde{M}_{h,2},M_{h,2})=codim(I_h\widetilde{M}_{h,2},I_hM_{h,2})\leqslant codim(\Delta_{h,2}^{-1}\widetilde{V}_{h0}, V_{h0}) \leqslant m_0$.

\begin{theorem}\label{thm:sdcfs}
We have a stable decomposition of $M_{h,k}$:
\begin{equation}\label{sdopefs}
\|\nabla^2_hw_h\|_{0,\Omega}^2\lesssim \sum_{T\in\mathcal{T}_h} h_T^{-4}\|w_h-\Pi_{h,k}I_hw_h\|_{0,T}^2+(\Delta_{h,k}I_hw_h,\Delta_{h,k}I_hw_h)\lesssim \|\nabla^2_hw_h\|_{0,\Omega}^2,
\end{equation}
which holds for all $w_h\in M_{h,k}$ and $k=1,2$, with this exception: when $k=2$,  the left inequality above only holds for $w_h\in \widetilde{M}_{h,2}$.
\end{theorem}
\begin{proof}
The theorem follows from Lemma \ref{lem:dlevc} and Lemma \ref{lem:sds} below.
\end{proof}

\begin{lemma}\label{lem:avesqr}
Let $m\geqslant1$ be an integer. The following equivalence for $\gamma,\beta_1,\dots,\beta_m\in\mathbb{R}$ depends on $m$ only:
$$
\frac{1}{m}\sum_{i=1}^m(\gamma-\beta_i)^2\cequiv (\gamma-\bar{\beta})^2+\frac{1}{m}\sum_{i=1}^m(\beta_i-\beta_{i+1})^2,
$$
where $\bar{\beta}=\frac{1}{m}(\beta_1+\dots+\beta_m)$, $\beta_{m+1}=\beta_1$.
\end{lemma}
The proof of Lemma \ref{lem:avesqr} is straightforward; therefore it is omitted.

\begin{lemma}\label{lem:linterr}
\cite{WangXu2006}
It holds for all $w_h\in M_{h,2}$ that
\begin{equation}
\|w_h-I_hw_h\|_{0,T}^2\lesssim \sum_{T'\subset\omega_T}h_{T'}^4|w_h|_{2,T'}^2,\ \ \forall\,T\in\mathcal{T}_h.
\end{equation}
\end{lemma}

\begin{lemma}\label{lem:sinterr}
For all $p_h\in V_{h0}$, it holds for $k=1,2$ that
\begin{equation}\label{eq:interrdfs}
\sum_{T\in\mathcal{T}_h}h_T^{-4}\|\Pi_{h,k}p_h-p_h\|_{0,T}^2\lesssim \sum_{e\in\mathcal{E}_{h,k}}h_e^{-1}\int_e\Big\llbracket\frac{\partial p_h}{\partial\mathbf{n}_e}\Big\rrbracket^2.
\end{equation}
\end{lemma}
\begin{proof}
By Assumptions A1 and A3, we have
\begin{equation}\label{cbnpc}
\sum_{T\in\mathcal{T}_h}h_T^{-4}\|\Pi_{h,k}p_h-p_h\|_{0,T}^2\cequiv \sum_{\deg(\alpha)=1}\sum_{T\subset\omega_\alpha}\Big|\mathcal{N}_\alpha(p_h|_T)-\mathcal{N}_\alpha(\Pi_{h,k}p_h)\Big|^2.
\end{equation}
Let $\mathcal{N}_\alpha$ be a nodal parameter with $deg(\alpha)=1$, and we rewrite the element patch $\omega_\alpha$ as $\omega_\alpha=\{T^\alpha_j\}_{j=1}^{\#\omega_\alpha}$, such that $T^\alpha_{j}$ and $T^\alpha_{j+1}$ share a common edge, $j=1,\dots,\#\omega_\alpha-1$. Denote $T_{\#\omega_\alpha+1}=T_1$, then
\begin{multline}\label{interior}
\sum_{T\subset\omega_\alpha}\Big|\mathcal{N}_\alpha(p_h|_T)-\mathcal{N}_\alpha(\Pi_{h,k}p_h)\Big|^2\\ \cequiv\#\omega_\alpha\Big|\mathcal{N}_\alpha(\Pi_{h,k}p_h)-\frac{1}{\#\omega_\alpha}\sum_{T\subset\omega_\alpha}\mathcal{N}_\alpha(p_h|_T)\Big|^2+\sum_{j=1}^{\#\omega_\alpha}\Big(\mathcal{N}_\alpha(p_h|_{T^\alpha_j})-\mathcal{N}_\alpha(p_h|_{T^\alpha_{j+1}})\Big)^2,\ \ k=1,2.
\end{multline}

If $D_\alpha\not\subset\partial\Omega$, then by the definitions of $\Pi_{h,k}$,
\begin{equation}
\sum_{T\subset\omega_\alpha}\Big|\mathcal{N}_\alpha(p_h|_T)-\mathcal{N}_\alpha(\Pi_{h,k}p_h)\Big|^2 \cequiv\sum_{j=1}^{\#\omega_\alpha}\Big(\mathcal{N}_\alpha(p_h|_{T^\alpha_j})-\mathcal{N}_\alpha(p_h|_{T^\alpha_{j+1}})\Big)^2,
\end{equation}
and by the continuity of the piecewise linear function $p_h$, we obtain further that
\begin{equation}
\sum_{T\subset\omega_\alpha}\Big|\mathcal{N}_\alpha(p_h|_T)-\mathcal{N}_\alpha(\Pi_{h,k}p_h)\Big|^2\lesssim \sum_{\overline{e}\cap\overline{D_\alpha}\neq\emptyset}h_e^{-1}\int_e\Big\llbracket\frac{\partial p_h}{\partial\mathbf{n}_e}\Big\rrbracket^2, \ \ k=1,2.
\end{equation}

If $D_\alpha\subset\partial\Omega$, then for $\Pi_{h,1}$ we obtain that\\
\begin{equation}
\sum_{T\subset\omega_\alpha}\big|\mathcal{N}_\alpha(p_h|_T)-\mathcal{N}_\alpha(\Pi_{h,1}p_h)\big|^2=\sum_{T\subset\omega_\alpha}\big|\mathcal{N}_\alpha(p_h|_T)\big|^2\lesssim \sum_{\overline{e}\cap\overline{D_\alpha}\neq\emptyset} h_e^{-1}\int_e\Big\llbracket\frac{\partial p_h}{\partial\mathbf{n}_e}\Big\rrbracket^2.
\end{equation}
For $\Pi_{h,2}$, there are two cases. If $deg(\alpha)=1$, then noting that $p_h|_\Gamma=0$, we obtain
\begin{equation}
\begin{array}{rl}
\displaystyle \sum_{T\in\omega_\alpha}\Big|\mathcal{N}_\alpha(p_h|_T)&\displaystyle-\mathcal{N}_\alpha(\Pi_{h,2}p_h)\Big|^2 \\
 = & \displaystyle \sum_{T\subset\omega_\alpha}\Big|\fint_{D_\alpha}\partial_\mathbf{n}(p_h|_T)(\mathbf{n}\cdot \mathbf{t}_\alpha)+\fint_{D_\alpha}\partial_\tau(p_h|_T)(\tau\cdot \mathbf{t}_\alpha)-\mathcal{N}_\alpha(\Pi_{h,2}p_h)\Big|^2\\
\lesssim & \displaystyle \sum_{T\subset\omega_\alpha}\Big|\fint_{D_\alpha}\partial_\mathbf{n}(p_h|_T)(\mathbf{n}\cdot \mathbf{t}_\alpha)-\mathcal{N}_\alpha(\Pi_{h,2}p_h)\Big|^2+\sum_{T\subset\omega_\alpha}\Big|\fint_{D_\alpha}\partial_\tau(p_h|_T)(\tau\cdot \mathbf{t}_\alpha)\Big|^2 \\
\cequiv & \displaystyle \sum_{j=1}^{\#\omega_\alpha}(\fint_{D_\alpha}\partial_\mathbf{n}(p_h|_{T_j})-\fint_{D_\alpha}\partial_\mathbf{n}(p_h|_{T_{j+1}}))^2(\mathbf{n}\cdot\mathbf{t}_\alpha)^2+\sum_{T\subset\omega_\alpha}\Big|\fint_{D_\alpha}\partial_\tau(p_h|_T)(\tau\cdot\mathbf{t}_\alpha)\Big|^2 \\
\lesssim & \displaystyle \sum_{\overline{e}\cap\overline{D_\alpha}\neq\emptyset,e\not\subset\partial\Omega}h_e^{-1}\int_e\Big\llbracket\frac{\partial p_h}{\partial\mathbf{n}_e}\Big\rrbracket^2.
\end{array}
\end{equation}
Otherwise, by definition, we have
\begin{equation}
\sum_{T\subset\omega_\alpha}\Big|\mathcal{N}_\alpha(p_h|_T)-\mathcal{N}_\alpha(\Pi_{h,2}p_h)\Big|^2 \cequiv\sum_{j=1}^{\#\omega_\alpha}\Big(\mathcal{N}_\alpha(p_h|_{T^\alpha_j})-\mathcal{N}_\alpha(p_h|_{T^\alpha_{j+1}})\Big)^2.
\end{equation}
And, by the continuity of $p_h$ again, we obtain that
\begin{equation}
\sum_{T\subset\omega_\alpha}\Big|\mathcal{N}_\alpha(p_h|_T)-\mathcal{N}_\alpha(\Pi_{h,2}p_h)\Big|^2\lesssim \sum_{\overline{e}\cap\overline{D_\alpha}\neq\emptyset,e\not\subset\partial\Omega}h_e^{-1}\int_e\Big\llbracket\frac{\partial p_h}{\partial\mathbf{n}_e}\Big\rrbracket^2.
\end{equation}

Summing all the inequalities above leads to \eqref{eq:interrdfs} for $k=1,2$.
\end{proof}

\begin{lemma}\label{lem:sds}
It holds for $k=1,2$ that
\begin{equation}\label{eqsdfundfs}
\|\nabla^2_hw_h\|_{0,\Omega}^2\cequiv \sum_{T\in\mathcal{T}_h}h_T^{-4}\|w_h-\Pi_{h,k}I_hw_h\|_{0,T}^2+\sum_{e\in\mathcal{E}_{h,k}}h_e^{-1}\int_e\Big\llbracket\frac{\partial I_hw_h}{\partial\mathbf{n}_e}\Big\rrbracket^2,\ \mbox{for}\ w_h\in\,M_{h,k}.
\end{equation}
\end{lemma}
\begin{proof}
Denote $w_h^I:=I_hw_h$ for any $w_h\in M_{h,k}$, $k=1,2$. By inverse inequality, we have that
$$
\|\nabla^2_hw_h\|_{0,\Omega}^2=\sum_{T\in\mathcal{T}_h}|w_h|_{2,T}^2=\sum_{T\in\mathcal{T}_h}|w_h-w_h^I|_{2,T}^2\lesssim \sum_{T\in\mathcal{T}_h}h_T^{-4}\|w_h-w_h^I\|_{0,T}^2.
$$
Thus, by Lemma \ref{lem:linterr} and Assumption A3,
\begin{equation}\label{16}
\begin{array}{rl}
|w_h|_{2,h}^2\cequiv &\displaystyle\sum_{T\in\mathcal{T}_h} h_T^{-4}\|w_h-w_h^I\|_{0,T}^2\\ \cequiv& \displaystyle\sum_{T\in\mathcal{T}_h}h_T^{-4}\sum_{\omega_\alpha\supset T}h_T^{2+2deg(\alpha)}\Big(\mathcal{N}_\alpha(w_h)-\mathcal{N}_\alpha(w_h^I|_T)\Big)^2 \\
\cequiv &\displaystyle\sum_\alpha h_\alpha^{2deg(\alpha)-2}\sum_{T\subset\omega_\alpha}\Big(\mathcal{N}_\alpha(w_h)-\mathcal{N}_\alpha(w_h^I|_T)\Big)^2.
\end{array}
\end{equation}
It is straightforward to obtain that, for any $\alpha$, any $w_h\in M_{h,k}$ and any $p_h\in V_{h0}$, $k=1,2$,
\begin{multline}
\sum_{T\subset\omega_\alpha}\Big(\mathcal{N}_\alpha(w_h)-\mathcal{N}_\alpha(p_h|_T)\Big)^2= \#\omega_\alpha\Big(\mathcal{N}_\alpha(w_h)-\mathcal{N}_\alpha(\Pi_{h,k}p_h)\Big)^2+\sum_{T\in\omega_\alpha}\Big(\mathcal{N}_\alpha(\Pi_{h,k}p_h)-\mathcal{N}_\alpha(p_h)\Big)^2.
\end{multline}
Then by Assumption A3 again, given any $w_h\in M_{h,k}$,
\begin{equation}
\begin{array}{rl}
|w_h|_{2,h}^2\cequiv &\displaystyle\sum_\alpha h_\alpha^{2deg(\alpha)-2}\sum_{T\subset\omega_\alpha}\Big[\Big(\mathcal{N}_\alpha(w_h)-\mathcal{N}_\alpha(\Pi_{h,k}w_h^I)\Big)^2\\ \qquad&\displaystyle+\Big(\mathcal{N}_\alpha(\Pi_{h,k}w_h^I)-\mathcal{N}_\alpha(w_h^I|_T)\Big)^2\Big]\\
\cequiv &\displaystyle\sum_{T\in\mathcal{T}_h} h_T^{-4}\|w_h-\Pi_{h,k}w_h^I\|_{0,T}^2+\sum_{T\in\mathcal{T}_h} h_T^{-4}\|\Pi_{h,k}w_h^I-w_h^I\|_{0,T}^2,
\end{array}
\end{equation}
therefore, by \eqref{eq:interrdfs},
$$
|w_h|_{2,h}^2\lesssim \sum_{T\in\mathcal{T}_h} h_T^{-4}\|w_h-\Pi_{h,k}w_h^I\|_{0,T}^2+\sum_{e\in\mathcal{E}_{h,k}}h_e^{-1}\int_e\Big\llbracket\frac{\partial w_h^I}{\partial\mathbf{n}_e}\Big\rrbracket^2\ \mbox{for}\,w_h\in M_{h,k},\ k=1,2.
$$
On the other hand, as
$$
\sum_{e\in\mathcal{E}_{h,k}}h_e^{-1}\int_e\Big\llbracket\frac{\partial w_h^I}{\partial\mathbf{n}_e}\Big\rrbracket^2\lesssim \sum_{e\in\mathcal{E}_{h,k}}h_e^{-1}\int_e\Big\llbracket\frac{\partial(w_h^I-w_h)}{\partial\mathbf{n}_e}\Big\rrbracket^2+\sum_{e\in\mathcal{E}_{h,k}}h_e^{-1}\int_e\Big\llbracket\frac{\partial w_h}{\partial\mathbf{n}_e}\Big\rrbracket^2\lesssim |w_h|_{2,h}^2,
$$
we obtain
$$
|w_h|_{2,h}^2\gtrsim \sum_{T\in\mathcal{T}_h} h_T^{-4}\|w_h-\Pi_{h,k}w_h^I\|_{0,T}^2+\sum_{e\in\mathcal{E}_{h,k}}h_e^{-1}\int_e\Big\llbracket\frac{\partial w_h^I}{\partial\mathbf{n}_e}\Big\rrbracket^2,\ k=1,2.
$$
Combining these two points, we obtain \eqref{eqsdfundfs}.
\end{proof}

\section{Trace spaces on the boundary and their discretizations}\label{sec:reg}
\subsection{Trace spaces related to $H^1(\Omega)$}

As the trace space of $H^1(\Omega)$, $H^{1/2}(\Gamma)$ is a Hilbert space with respect to the norm $\displaystyle\|\lambda\|_{1/2,\Gamma}=\inf_{w\in H^1(\Omega),\lambda=w|_\Gamma}\|w\|_{1,\Omega}$. Denote $H^{-1/2}(\Gamma)$ as the dual space of $H^{1/2}(\Gamma)$.

Given $\lambda\in H^{1/2}\big(\Gamma\big)$, there exists a unique  $u_\lambda\in H^1(\Omega)$ satisfying
$$
\left\{
\begin{array}{rl}
-\Delta u_\lambda=0&\mbox{in}\,\Omega\\
u_\lambda=\lambda&\mbox{on}\,\partial\Omega.
\end{array}
\right.
$$
This defines a harmonic extension operator by $E\lambda=u_\lambda$. Then, $\|E\lambda\|_{1,\Omega}\cequiv \|\lambda\|_{1/2,\Gamma}$\  for $\lambda\in H^{1/2}(\Gamma)$.

For any edge $\Gamma_i$ of $\Omega$, define the space $H^{1/2}_{00}(\Gamma_i):=\{\lambda\in L^2(\Gamma_i):\tilde{\lambda}\in H^{1/2}(\Gamma)\}$ where $\tilde{\lambda}$ is the zero extension of $\lambda$ into $\Gamma\setminus\Gamma_i$, and the space is a Hilbert space with a norm given by $\|\lambda\|_{H^{1/2}_{00}(\Gamma_i)}:=\|\tilde{\lambda}\|_{1/2,\Gamma}$. Moreover, $H^{1/2}_{00}(\Gamma_i)$ is the interpolated space halfway between the $H^1_0(\Gamma_i)$ and $L^2(\Gamma_i)$ spaces. We refer to \cite{LionsMagenes1972} and references therein for more details.

\subsection{Trace spaces related to $H^2(\Omega)\cap H^1_0(\Omega)$}

We consider the following subspace of $H^{1/2}(\Gamma)$:
\begin{equation}
H^{1/2}_c(\Gamma):=\{\lambda\in H^{1/2}(\Gamma):\lambda\chi_i\in H^{1/2}(\Gamma),\ i=1,2,\dots,K\},
\end{equation}
where $\chi_i$ is the characteristic function on $\Gamma_i$, with the norm
$$
\|\lambda\|_{1/2,c,\Gamma}:=(\sum_{i=1}^K\|\lambda\chi_i\|_{1/2,\Gamma}^2)^{1/2}.
$$
\begin{proposition}
For any $\lambda\in L^2(\Gamma)$, the sufficient and necessary condition of $\lambda\in H^{1/2}_c(\Gamma)$ is $\lambda|_{\Gamma_i}\in H^{1/2}_{00}(\Gamma_i)$ for $i=1,\dots,K$. Moreover, $\|\lambda\|_{1/2,c,\Gamma}\cequiv\sum_{i=1}^K\|\lambda|_{\Gamma_i}\|_{H^{1/2}_{00}(\Gamma_i)}$, and $H^{1/2}_c(\Gamma)$ is a Hilbert space with respect to $\|\cdot\|_{1/2,c,\Gamma}$.
\end{proposition}

\begin{theorem}
The following identities hold both algebraically and topologically:
\begin{equation}\label{eq:algdes}
\begin{array}{rl}
\displaystyle H^{1/2}_c(\Gamma)&=\displaystyle \Bigg\{\lambda\in L^2(\Gamma): \lambda\mathbf{n}\in (H^{1/2}(\Gamma))^2\Bigg\}\displaystyle=\Bigg\{\lambda\in L^2(\Gamma): \lambda\mathbf{\tau}\in (H^{1/2}(\Gamma))^2\Bigg\}.
\end{array}
\end{equation}
More specifically, $\|\lambda\|_{1/2,c,\Gamma}\cequiv\|\lambda\mathbf{n}\|_{1/2,\Gamma}=\|\lambda\mathbf{\tau}\|_{1/2,\Gamma}$ for $\lambda\in H^{1/2}_c(\Gamma)$.
\end{theorem}

\begin{proof}
To begin, we prove the first identity in \eqref{eq:algdes}. Let $\lambda\in H_c^{1/2}(\Gamma)$.
We have $\lambda\mathbf{n}=\sum_{i=1}^K\lambda\chi_i\mathbf{n}_{\Gamma_i}$. By definition, $\lambda\chi_i\in H^{1/2}(\Gamma)$; hence,
$\lambda\chi_i\mathbf{n}_{\Gamma_i}\in (H^{1/2}(\Gamma))^2$; thus,
$\lambda\mathbf{n}\in  (H^{1/2}(\Gamma))^2$.

Now let $\lambda\in L^2(\Gamma)$ be such that
$\lambda\mathbf{n}\in  (H^{1/2}(\Gamma))^2$.   It can be easily verified for $j=i-1$ and $i+1$ that $\lambda \chi_i=\lambda\mathbf{n}\cdot \frac{\tau_{\Gamma_j}}{\mathbf{n}_{\Gamma_i}\cdot\tau_{\Gamma_j}}$ on $\Gamma_i\cup\Gamma_j$ and $\lambda\chi_i=\lambda\mathbf{n}\cdot\frac{\tau_{\Gamma_j}}{\mathbf{n}_{\Gamma_i}\cdot\tau_{\Gamma_j}}=0 \ \mbox{on\ } \Gamma_j$. Namely, $\lambda\chi_i$ can be expressed as
a linear combination of function $\mathbf{0}$ and functions $\lambda\mathbf{n}\cdot \frac{\tau_{\Gamma_j}}{\mathbf{n}_{\Gamma_i}\cdot\tau_{\Gamma_j}}$ for
$j=i\pm1$, which belong to $H^{1/2}(\Gamma)$, with the coefficients in $C^{\infty}(\Gamma)$. Therefore, $\lambda\chi_i$ itself belongs to $H^{1/2}(\Gamma)$.

It is easy to verify that $\|\lambda\mathbf{n}\|_{1/2,\Gamma}$ is a norm on
$\lambda\in H^{1/2}_c(\Gamma)$. Moreover,
$\|\lambda\mathbf{n}\|_{1/2,\Gamma}\lesssim
\|\lambda\|_{1/2,c,\Gamma}$. It is straightforward to verify that
$H^{1/2}_c(\Gamma)$ is a Hilbert space with respect to both the norms
$\|\cdot\mathbf{n}\|_{1/2,\Gamma}$ and $\|\cdot\|_{1/2,c,\Gamma}$,
therefore, by open mapping theorem, $\|\lambda\|_{1/2,c,\Gamma}\cequiv \|\lambda\mathbf{n}\|_{1/2,\Gamma}$ for $\lambda\in H^{1/2}_c(\Gamma)$. This finishes the proof of the first identity.

The second identity in \eqref{eq:algdes} is obvious given that $\tau=\mathbf{n}^\perp$. And, the theorem is thus proved.
\end{proof}

Denote $H^{-1/2}_c\big(\Gamma\big)$ as the dual space of $H^{1/2}_c\big(\Gamma\big)$, with the norm $\displaystyle\|\chi\|_{-1/2,c,\Gamma}=\sup_{\lambda\in H^{1/2}_c\big(\Gamma\big)}\frac{\langle\chi,\lambda\rangle}{\|\lambda\|_{1/2,c,\Gamma}}$, where $\langle\cdot,\cdot\rangle$ is the duality between $H^{-1/2}_c\big(\Gamma\big)$ and  $H^{1/2}_c\big(\Gamma\big)$.

Define a space of biharmonic functions $\widetilde{H}^2(\Omega):=\{u\in H^2(\Omega)\cap H^1_0(\Omega): (\Delta u,\Delta v)=0,\ \forall\,v\in H^2_0(\Omega)\}$. It can be verified that $\nabla(\widetilde{H}^2(\Omega))$ is a Hilbert subspace of $H(div;\Omega)$ and $\Delta(\widetilde{H}^2(\Omega))$ is a Hilbert subspace of $L^2(\Omega)$. In the remainder of this section, we show that $H^{1/2}_c(\Gamma)$ and $H^{-1/2}_c(\Gamma)$ are isomorphic trace spaces of $\nabla(\widetilde{H}^2(\Omega))$ and $\Delta(\widetilde{H}^2(\Omega))$, respectively.

\begin{theorem}\label{lem:trace}
$H^{1/2}_c(\Gamma)$ is the normal derivative trace space of $H^2(\Omega)\cap H^1_0(\Omega)$ in the sense that:
\begin{enumerate}
\item If $u\in H^2(\Omega)\cap H^1_0(\Omega)$, then $\frac{\partial u}{\partial\mathbf{n}}\big|_\Gamma\in H^{1/2}_c\big(\Gamma\big)$ and $\|\frac{\partial u}{\partial\mathbf{n}}\|_{1/2,c,\Gamma}\lesssim\|u\|_{2,\Omega}$.
\item Given any $\lambda\in H^{1/2}_c\big(\Gamma\big)$, there exists a unique $u\in \widetilde{H}^2(\Omega)$, such that $\frac{\partial u}{\partial\mathbf{n}}|_\Gamma=\lambda$ and $\|u\|_{2,\Omega}\lesssim\|\lambda\|_{1/2,c,\Gamma}$.
\end{enumerate}
\end{theorem}
\begin{proof}
Given $u\in H^2(\Omega)\cap H^1_0(\Omega)$, then $\nabla u\in \Big(H^1\big(\Omega\big)\Big)^2$ and $\nabla u|_\Gamma\in \Big(H^{1/2}\big(\Gamma\big)\Big)^2$. As $\frac{\partial u}{\partial \tau}=0$ along $\Gamma$, it holds that $\frac{\partial u}{\partial\mathbf{n}}\mathbf{n}=\nabla u$ and $\frac{\partial u}{\partial\mathbf{n}}\in H^{1/2}_c\big(\Gamma\big)$. Moreover, $\big\|\frac{\partial u}{\partial\mathbf{n}}\big\|_{1/2,c,\Gamma}=\|\nabla u\|_{1/2,\Gamma}\lesssim \|u\|_{2,\Omega}$. This proves the first part of the theorem.

To prove the second part, we consider an auxiliary Stokes problem for any $\lambda\in H^{1/2}_c\big(\Gamma\big)$:
\begin{equation*}
\left\{\begin{array}{rl}
\Delta\mathbf{\psi}+\nabla p=0&\mbox{in}\,\Omega,\\
\nabla\cdot\mathbf{\psi}=0&\mbox{in}\,\Omega,\\
\mathbf{\psi}=\lambda\tau&\mbox{on}\,\Gamma.
\end{array}
\right.
\end{equation*}
This problem, thanks to the trivial fact of $\tau\cdot\mathbf{n}=0$ on $\Gamma$, obviously admits a unique solution, such that $(\mathbf{\psi},p)\in \Big(H^1\big(\Omega\big)\Big)^2\times L^2_0(\Omega)$ satisfies $\|\mathbf{\psi}\|_{1,\Omega}\lesssim \|\lambda\tau\|_{1/2,\Gamma}$. (see Theorem 5.1 and Remark 5.3 on pp. 80-83 in \cite{GiraultRaviart1986}). As $\nabla\cdot\mathbf{\psi}=0$, there exists a $u\in H^1(\Omega)$ such that $\psi=\mathrm{curl} u$ (see Theorem 3.1 on p. 37 in \cite{GiraultRaviart1986}). And,  $\frac{\partial u}{\partial\mathbf{n}}=\mathrm{curl}u\cdot\tau=\mathbf{\psi}\cdot\mathsf{\tau}=\lambda$ on $\Gamma$. Given $v\in H^2_0(\Omega)$, $(\Delta u,\Delta v)=(\nabla\times\mathrm{curl}u,\nabla\times\mathrm{curl}v)=(\nabla\times\psi,\nabla\times\mathrm{curl}v)=(\nabla\times\psi,\nabla\times\mathrm{curl}v)+(p,\nabla\cdot\mathrm{curl}v)=0$. Further, $\frac{\partial u}{\partial \tau}=\nabla u\cdot\tau=\mathbf{\psi}\cdot\mathbf{n}=0$ along $\Gamma$; therefore, $u$ is a constant along $\Gamma$, and thus we may choose $u\in H^1_0(\Omega)$. Then, $u\in \widetilde{H}^2(\Omega)$ and $\|u\|_{2,\Omega}\lesssim|u|_{1,\Omega}+|u|_{2,\Omega}\cequiv \|\nabla u\|_{1,\Omega}=\|\mathbf{\psi}\|_{1,\Omega}\lesssim \|\lambda\tau\|_{1/2,\Gamma}\cequiv\|\lambda\|_{1/2,c,\Gamma}$.  The uniqueness of such a $u$ is straightforward. This finishes the proof of the theorem.
\end{proof}

By Theorem \ref{lem:trace}, given $\lambda\in H^{1/2}_c(\Gamma)$, there exists a unique $u_\lambda\in \widetilde{H}^2(\Omega)$ such that $\frac{\partial u_\lambda}{\partial\mathbf{n}}=\lambda$. This defines an extension operator $E_d:H^{1/2}_c(\Gamma)\to\nabla(\widetilde{H}^2(\Omega))$ by $E_d\lambda:=\nabla u_\lambda$. Moreover, we can define a trace operator $T_r^d:\nabla (\widetilde{H}^2(\Omega))\to H^{1/2}_c(\Gamma)$ by $T_r^d(\nabla u)=\nabla u\cdot\mathbf{n}$ for $u\in\widetilde{H}^2(\Omega)$. The following lemma follows from Theorem \ref{lem:trace} directly.

\begin{lemma}\label{lem:isomphb}
 $T_r^d$ is an isomorphism from $(\nabla(\widetilde{H}^2(\Omega)),\|\cdot\|_{div})$ onto $H^{1/2}_c(\Gamma)$, $E_d$ is an isomorphism from $H^{1/2}_c(\Gamma)$ onto $(\nabla\widetilde{H}^2(\Omega),\|\cdot\|_{div})$, and $T_r^d\circ E_d=Id_{H^{1/2}_c(\Gamma)}$ and $E_d\circ T_r^d=Id_{\nabla (\widetilde{H}^2(\Omega))}$.
\end{lemma}

\begin{remark}
It is easy to verify that for $u\in H^2(\Omega)\cap H^1_0(\Omega)$, $|\nabla u|_{div}=|\nabla u|_{1,\Omega}$, therefore, $\|\cdot\|_{1,\Omega}$ and $\|\cdot\|_{div}$ are equivalent on $\nabla(H^2(\Omega)\cap H^1_0(\Omega))$. Let $\widetilde{T}_r^d:\nabla (H^2(\Omega)\cap H^1_0(\Omega))\to H^{1/2}_c(\Gamma)$ be defined by $\widetilde{T}_r^d(\nabla u)=\nabla u\cdot\mathbf{n}$ for $u\in H^2(\Omega)\cap H^1_0(\Omega)$. Then the range of $\widetilde{T}_r^d$ is $H^{1/2}_c(\Gamma)$, and $\widetilde{T}_r^d$ is a continuous extension of $T_r^d$.
\end{remark}

\begin{theorem}\label{lem:bihext}
The spaces $H^{-1/2}_c(\Gamma)$ and $\Delta(\widetilde{H}^2(\Omega))$ are isomorphic in the sense that
\begin{enumerate}
\item For any $w\in\widetilde{H}^2(\Omega)$, there exists a unique $\zeta_w\in H^{-1/2}_c(\Gamma)$, such that
\begin{equation}
\langle\zeta_w,\frac{\partial v}{\partial\mathbf{n}}\rangle=(\Delta w,\Delta v),\quad\forall\,v\in H^2(\Omega)\cap H^1_0(\Omega),
\end{equation}
and $\|\zeta_w\|_{-1/2,c,\Gamma}\lesssim \|\Delta w\|_{0,\Omega}$.
\item For any $\zeta\in H^{-1/2}_c\big(\Gamma\big)$, there exists a unique $w_\zeta\in \widetilde{H}^2(\Omega)$, such that
\begin{equation}\label{eq:dualext}
(\Delta w_\zeta,\Delta v)=\langle\zeta,\frac{\partial v}{\partial\mathbf{n}}\rangle,\quad\forall\, v\in H^2(\Omega)\cap H^1_0(\Omega),
\end{equation}
and $\|\Delta w_\zeta\|_{0,\Omega}\lesssim \|\zeta\|_{-1/2,c,\Gamma}$.
\end{enumerate}
\end{theorem}

\begin{proof}
Let us first prove the first part of the theorem. Given $\lambda\in H^{1/2}_c(\Gamma)$, then by Theorem \ref{lem:trace}, there exists a $v_\lambda\in\widetilde{H}^2(\Omega)$, such that $\lambda=\frac{\partial v_\lambda}{\partial\mathbf{n}}|_\Gamma$. We define a linear functional $\zeta_w$ on $H^{1/2}_c(\Gamma)$ by $\langle\zeta_w,\lambda\rangle=(\Delta w,\Delta v_\lambda)$. For any $\mu\in H^{1/2}_c\big(\Gamma\big)$, denote $v_\mu\in H^2(\Omega)\cap H^1_0(\Omega)$ such that $\frac{\partial v_\mu}{\partial\mathbf{n}}|_\Gamma=\mu$ and $\|v_\mu\|_{2,\Omega}\leqslant C\|\mu\|_{1/2,c,\Gamma}$ with a generic constant $C$. Then,
$$
\|\zeta\|_{-1/2,c,\Gamma}=\sup_{\mu\in H^{1/2}_c\big(\Gamma\big)}\frac{\langle\zeta,\mu\rangle}{\|\mu\|_{1/2,c,\Gamma}}=\sup_{\mu\in H^{1/2}_c\big(\Gamma\big)}\frac{(\Delta w_\zeta,\Delta v_\mu)}{\|\mu\|_{1/2,c,\Gamma}}\leqslant C \|\Delta w_\zeta\|_{0,\Omega}.
$$

To prove the second part of the theorem, let $\zeta\in H^{-1/2}_c\big(\Gamma\big)$. Then, there is a unique $w_\zeta\in\widetilde{H}^2(\Omega) $ such that $(\Delta w_\zeta,\Delta v)=\langle\zeta,\frac{\partial v}{\partial\mathbf{n}}\rangle$, $\forall\,v\in H^2(\Omega)\cap H^1_0(\Omega)$. Furthermore, $(\Delta w_\zeta,\Delta w_\zeta)=\langle\zeta,\frac{\partial w_\zeta}{\partial\mathbf{n}}\rangle\lesssim \|\zeta\|_{-1/2,c,\Gamma}\|\frac{\partial w_\zeta}{\partial\mathbf{n}}\|_{1/2,c,\Gamma}\lesssim \|\zeta\|_{-1/2,c,\Gamma}\|w_\zeta\|_{2,\Omega}\lesssim \|\zeta\|_{-1/2,c,\Gamma}\|\Delta w_\zeta\|_{0,\Omega}$; therefore, $\|\Delta w_\zeta\|_{0,\Omega}\lesssim \|\zeta\|_{-1/2,c,\Gamma}$. The proof is finished.
\end{proof}

By Theorem \ref{lem:bihext}, given $\zeta\in H^{-1/2}_c(\Gamma)$, there exists a unique $w_\zeta\in \widetilde{H}^2(\Omega)$ such that $(\Delta w_\zeta,\Delta v)=\langle\zeta,\frac{\partial v}{\partial\mathbf{n}}\rangle$ for all $v\in H^2(\Omega)\cap H^1_0(\Omega)$. We define an extension operator $E_c:H^{-1/2}_c(\Gamma)\to \Delta(\widetilde{H}^2(\Omega))$ by $E_c\zeta:=\Delta w_\zeta$. Moerover, we define a trace operator $T_r^c:\Delta(\widetilde{H}^2(\Omega))\to H^{-1/2}_c(\Gamma)$ by $\langle T_r^c (\Delta w),\frac{\partial v}{\partial\mathbf{n}}\rangle=(\Delta w,\Delta v)\ \forall\,v\in H^2(\Omega)\cap H^1_0(\Omega)$ for $w\in\widetilde{H}^2(\Omega)$. The following lemma follows from Theorem \ref{lem:bihext} directly.

\begin{lemma}\label{lem:isomnhb}
 $T_r^c$ is an isomorphism from $(\Delta(\widetilde{H}^2(\Omega)),\|\cdot\|_{0,\Omega})$ onto $H^{-1/2}_c(\Gamma)$, $E_c$ is an isomorphism from $H^{-1/2}_c(\Gamma)$ onto $(\Delta(\widetilde{H}^2(\Omega)),\|\cdot\|_{0,\Omega})$, and $T_r^c\circ E_c=Id_{H^{-1/2}_c(\Gamma)}$ and $E_c\circ T_r^c=Id_{\Delta (\widetilde{H}^2(\Omega))}$.
\end{lemma}

It is easy to verify that if $\lambda\in H^{1/2}(\Gamma)$, then $(E\lambda,\Delta v)=(E_c\lambda,\Delta v),\ \forall\,v\in H^2(\Omega)\cap H^1_0(\Omega)$; that is, $E_c\lambda$ is indeed the $L^2$-projection of $E\lambda$ in $\mathcal{R}(\Delta)$. In particular, when $\mathcal{R}(\Delta)=L^2(\Omega)$, which holds when $\Omega$ is convex, $E_c$ and $E$ coincide on $H^{1/2}(\Gamma)$. Therefore, $E_c:H^{-1/2}_c(\Gamma)\to\mathcal{R}(\Delta)$ is a generalization of the harmonic extension operator to a larger function space $H^{-1/2}_c(\Gamma)\supset H^{1/2}(\Gamma)$, and $E_c$ is called a generalized harmonic extension operator.
\subsection{Generalized Poincar\'e--Steklov operator and its inverse}

Given $\lambda\in H^{1/2}_c(\Gamma)$, let $w_\lambda\in \widetilde{H}^2(\Omega)$ be such that $\nabla w_\lambda=E_d\lambda$. Define a generalized Poincar\'e--Steklov operator as $T_c:H^{1/2}_c(\Gamma)\to H^{-1/2}_c(\Gamma)$ by $T_c\lambda:=T_r^c(\Delta w_\lambda)$. Then
\begin{equation}\label{gcgps}
\langle T_c\lambda,\mu\rangle=(\nabla E_d\lambda,\nabla E_d\mu),\ \ \forall\,\lambda,\mu\in H^{1/2}_c(\Gamma).
\end{equation}

Given $\zeta\in H^{-1/2}_c(\Gamma)$, let $w_\zeta\in\widetilde{H}^2(\Omega)$ be such that $\Delta w_\zeta=E_c\zeta$. Define an operator $S_c$ on $H^{-1/2}_c(\Gamma)$ by $S_c\zeta:=T_r^d(\nabla w_\zeta)$. We directly obtain that
\begin{equation}\label{gcgips}
\langle\zeta,S_c\xi\rangle=(E_c\zeta,E_c\xi),\quad\forall\,\zeta,\xi\in H^{-1/2}_c(\Gamma).
\end{equation}
\begin{lemma}\label{lem:gips}
The operators $T_c$ and $S_c$ are both algebraic and topological isomorphisms between $H^{-1/2}_c(\Gamma)$ and $H^{1/2}_c(\Gamma)$, and $T_c\circ S_c=Id_{H^{-1/2}_c(\Gamma)}$ and $S_c\circ T_c=Id_{H^{1/2}_c(\Gamma)}$. Moreover, $\langle\zeta,S_c\zeta\rangle\cequiv \|\zeta\|_{-1/2,c,\Gamma}^2$ for $\zeta\in H^{-1/2}_c(\Gamma)$ and $\langle T_c\lambda,\lambda\rangle\cequiv \|\lambda\|_{1/2,c,\Gamma}^2$ for $\lambda\in H^{1/2}_c(\Gamma)$.
\end{lemma}
\begin{proof}
Firstly we show that $T_c$ is a bijection. It is straightforward to verify that $T_c$ is injective. Now let $\mu\in H^{-1/2}_c(\Gamma)$, then there exists a $w_\mu\in\widetilde{H}^2(\Omega)$, such that $\Delta w_\mu=E_c\mu$. Denote $\lambda=T_r^d\nabla w_\mu$, then by definition $\mu=T_c\lambda$. Thus $T_c$ is surjective. Similarly, $S_c$ is a bijection. By definition and the fact that $E_c$ and $T_r^c$ are inverse to each other and $E_d$ and $T_r^d$ are inverse to each other, $T_c\circ S_c=Id_{H^{-1/2}_c(\Gamma)}$ and $S_c\circ T_c=Id_{H^{1/2}_c(\Gamma)}$.

For any $\zeta\in H^{-1/2}_c(\Gamma)$, it can be derived by Theorems \ref{lem:trace} and \ref{lem:bihext} that $\|S_c\zeta\|_{1/2,c,\Gamma}\lesssim \|\zeta\|_{-1/2,c,\Gamma}$. Similarly, we can prove $\|T_c\lambda\|_{-1/2,c,\Gamma}\lesssim \|\lambda\|_{1/2,c,\Gamma}$ for $\lambda\in H^{1/2}_c(\Gamma)$.

The proof is finally finished by noting \eqref{gcgps} together with Lemma \ref{lem:isomnhb} and \eqref{gcgips} together with Lemma \ref{lem:isomphb}.
\end{proof}

By means of the inverse generalized Poincar\'e--Steklov operator $S_c$, the first biharmonic problem can be decomposed to two second biharmonic problems, as in the proposition below.
\begin{proposition}\label{lem:decouec}
Let $u\in H^2_0(\Omega)$ solve $(\Delta u,\Delta v)=(f,v)\ \forall\,v\in H^2_0(\Omega)$, then $u$ can be obtained by seeking $(\tilde{u},\zeta,u)\in (H^2(\Omega)\cap H^1_0(\Omega))\times H^{-1/2}_c(\Gamma)\times (H^2(\Omega)\cap H^1_0(\Omega)$), such that
\begin{eqnarray}
&(\Delta \tilde{u},\Delta p)=(f,p),&\forall\, p\in H^2(\Omega)\cap H^1_0(\Omega),\label{gpdis1}\\
&\langle\zeta,S_c\gamma\rangle=-(\Delta\tilde{u},E_c\gamma),&\forall\,\gamma\in H^{-1/2}_c(\Gamma),\label{gpdis2}\\
&(\Delta u,\Delta v)=(\Delta\tilde{u}+E_c\zeta,\Delta v),&\forall\,v\in H^2(\Omega)\cap H^1_0(\Omega).\label{gpdis3}
\end{eqnarray}
\end{proposition}
\begin{proof}
Let $\gamma\in H^{-1/2}_c(\Gamma)$, then by \eqref{gpdis3} and \eqref{gpdis2}, we have $(\Delta u,E_c\gamma)=(\Delta\tilde{u}+E_c\zeta,E_c\gamma)=0$. Namely, $\langle\gamma,\frac{\partial u}{\partial\mathbf{n}}\rangle=(E_c\gamma,\Delta u)=0$ for any $\gamma\in H^{-1/2}_c(\Gamma)$, thus $\frac{\partial u}{\partial\mathbf{n}}|_\Gamma=0$, and $u\in H^2_0(\Omega)$. Further, for any $v\in H^2_0(\Omega)$,
$$
(\Delta u,\Delta v)=(\Delta\tilde{u},\Delta v)+(E_c\zeta,\Delta v)=(f,v)+\langle\zeta,\frac{\partial v}{\partial\mathbf{n}}\rangle=(f,v).
$$
This finishes the proof.
\end{proof}

\subsection{Linear element spaces on the boundary}

Given a subset $F\subset\Omega$, define a restriction operator as $I_F^0:\mathcal{B}_h\to\mathcal{B}_h$ by $I_F^0v(x)=v(x)$ if $x$ is a vertex in $F$, and $I_F^0v(x)=0$ if $x$ is a vertex out of $F$. The lemma below is about the stability of $I_F^0$. 

\begin{lemma}\label{lem:XuZou}
\cite{XuZou1998} Let $F$ be a vertex on $\partial\Omega$ or an edge of $\partial\Omega$, then $\|I_F^0\lambda_h\|_{1/2,\Gamma}\lesssim (1+|\log h|)\|\lambda_h\|_{1/2,\Gamma}$.
\end{lemma}

\begin{lemma}\label{lem:-12}
It holds for $\lambda_h\in\mathcal{B}_h$ that
\begin{enumerate}
\item $\displaystyle\|\lambda_h\|_{-1/2,\Gamma}\cequiv \sup_{\mu_h\in\mathcal{B}_h}\frac{(\lambda_h,\mu_h)}{\|\mu_h\|_{1/2,\Gamma}}$;
\item $\displaystyle\|\lambda_h\|_{0,\Gamma}\lesssim h^{-1/2}\|\lambda_h\|_{-1/2,\Gamma}$;
\item $\displaystyle\|\lambda_h^c\|_{-1/2,\Gamma}\lesssim h^{1/2}(1+|\log h|)\|\lambda_h^c\|_{0,\Gamma}$.
\end{enumerate}
\end{lemma}
\begin{proof}
Let $\widetilde{\mathcal{Q}}_h$ be the $L^2$ projection to $\mathcal{B}_h$, then $\|\widetilde{\mathcal{Q}}_h\lambda\|_{1/2,\Gamma}\lesssim \|\lambda\|_{1/2,\Gamma}$ for $\lambda\in H^{1/2}(\Gamma)$. Therefore $\|\lambda_h\|_{-1/2,\Gamma}=\sup_{\mu\in H^{1/2}(\Gamma)}\frac{(\lambda_h,\mu)_\Gamma}{\|\mu\|_{1/2,\Gamma}}\lesssim \sup_{\mu\in H^{1/2}(\Gamma)}\frac{(\lambda_h,\widetilde{\mathcal{Q}}_h\mu)}{\|\widetilde{\mathcal{Q}}_h\mu\|}\lesssim\sup_{\mu_h\in\mathcal{B}_h}\frac{(\lambda_h,\mu_h)}{\|\mu_h\|_{1/2,\Gamma}}\leqslant \|\lambda_h\|_{-1/2,\Gamma}$. This proves the first item. The second item follows from the first item and the inverse inequality.

Now, any $\mu_h\in\mathcal{B}_h$ can be decomposed as $\mu_h=\mathring{\mu}_h+\mu_h^c$. As $\|\mu_h\|_{1/2,\Gamma}\geqslant\|\mathring{\mu}_h\|_{1/2,\Gamma}-\|\mu_h^c\|_{1/2,\Gamma}$ and $\|\mu_h^c\|_{1/2,\Gamma}\lesssim(1+|\log h|)\|\mu_h\|_{1/2,\Gamma}$, we obtain that $\|\mathring{\mu}_h\|_{1/2,\Gamma}+\|\mu_h^c\|_{1/2,\Gamma}\lesssim (1+|\log h|)\|\mu_h\|_{1/2,\Gamma}$. Then, as $\mathring{\mu}_h$ vanishes at the corner, $(\lambda_h^c,\mathring{\mu}_h)_\Gamma\lesssim\|\lambda_h^c\|_{0,\Gamma}(h\|\mathring{\mu}_h\|_{1,supp(\lambda_h^c)})\lesssim h^{1/2}\|\lambda_h^c\|_{0,\Gamma}\|\mathring{\mu}_h\|_{1/2,\Gamma}$, and $(\lambda_h^c,\mu_h^c)_\Gamma\lesssim \|\lambda_h^c\|_{0,\Gamma}(h\|\mu_h^c\|_{1,\Gamma})\lesssim h^{1/2}\|\lambda_h^c\|_{0,\Gamma}\|\mu_h^c\|_{1/2,\Gamma}$, we obtain that
$$
\|\lambda_h^c\|_{-1/2,\Gamma}\cequiv\sup_{\mu_h\in\mathcal{B}_h}\frac{(\lambda_h^c,\mu_h)_\Gamma}{\|\mu_h\|_{1/2,\Gamma}}\lesssim(1+|\log h|)\sup_{\mu_h\in\mathcal{B}_h}\frac{(\lambda_h^c,\mathring{\mu_h})_\Gamma+(\lambda_h^c,\mu_h^c)_\Gamma}{\|\mu_h^c\|_{1/2,\Gamma}+\|\mathring{\mu}_h\|_{1/2,\Gamma}}\lesssim h^{1/2}(1+|\log h|)\|\lambda_h^c\|_{0,\Gamma}.
$$
The last item is proved, and the proof of the lemma is finished.
\end{proof}

\begin{lemma}
Let $\Gamma_i$ be an edge of $\partial\Omega$. Denote $\mathring{\mathcal{B}}_{hi}:=\mathring{\mathcal{B}}_h|_{\Gamma_i}= \mathring{\mathcal{B}}_h\cap H^1_0(\Gamma_i)$, and define $\mathcal{Q}_{hi}$ the $L^2$ projection to $\mathring{\mathcal{B}}_{hi}$. Then,
\begin{equation}\label{eq:staqhi}
\|\lambda_i-\mathcal{Q}_{hi}\lambda_i\|_{0,\Gamma_i}+h^{1/2}\|\mathcal{Q}_{hi}\lambda_i\|_{H^{1/2}_{00}(\Gamma_i)}\lesssim h^{1/2}\|\lambda_i\|_{H^{1/2}_{00}(\Gamma)},\quad\forall\,\lambda_i\in H^{1/2}_{00}(\Gamma_i).
\end{equation}
\end{lemma}
\begin{proof}
Firstly, it can be proved that
\begin{equation}\label{eq:staql2}
\left\{
\begin{array}{ll}
\displaystyle\|\lambda_i-\mathcal{Q}_{hi}\lambda_i\|_{0,\Gamma_i}+\|\mathcal{Q}_{hi}\lambda_i\|_{0,\Gamma_i}\lesssim\|\lambda_i\|_{0,\Gamma_i},\quad\forall\,\lambda_i\in L^2(\Gamma_i);&\\
\displaystyle\|\lambda_i-\mathcal{Q}_{hi}\lambda_i\|_{0,\Gamma_i}+h\|\mathcal{Q}_{hi}\lambda_i\|_{1,\Gamma_i}\lesssim h\|\lambda_i\|_{1,\Gamma_i},\quad\forall\,\lambda_i\in H^1_{0}(\Gamma_i).&
\end{array}
\right.
\end{equation}
Then, as $H^{1/2}_{00}(\Gamma_i)$ is the interpolated space halfway between the $H^1_0(\Gamma_i)$ and $L^2(\Gamma_i)$ spaces, the result is obtained by interpolation.
\end{proof}

\begin{lemma}
Define $\mathcal{Q}_h$ as the $L^2$ projection to $\mathring{\mathcal{B}}_h$. Then
\begin{equation}\label{eq:staqh}
\|\lambda-\mathcal{Q}_h\lambda\|_{0,\Gamma}+h^{1/2}\|\mathcal{Q}_h\lambda\|_{1/2,c,\Gamma}\lesssim h^{1/2}\|\lambda\|_{1/2,c,\Gamma},\quad\forall\,\lambda\in H^{1/2}_c(\Gamma).
\end{equation}
\end{lemma}
\begin{proof}
For any $\lambda\in H^{1/2}_c(\Gamma)$, by definition, $(\mathcal{Q}_h\lambda)|_{\Gamma}=\mathcal{Q}_{hi}(\lambda|_{\Gamma_i})$. Denote $\lambda_i=\lambda|_{\Gamma_i}$ for $i=1,\dots,K$; then, $\|\lambda\|_{1/2,c,\Gamma} \cequiv \sum_{i=1}^K\|\lambda_i\|_{H^{1/2}_{00}(\Gamma_i)}$, $\|\mathcal{Q}_h\lambda\|_{1/2,c,\Gamma}\cequiv\sum_{i=1}^K\|\mathcal{Q}_{hi}\lambda_i\|_{H^{1/2}_{00}(\Gamma_i)}$ and $\|\lambda-\mathcal{Q}_h\lambda\|_{0,\Gamma}\cequiv \sum_{i=1}^K\|\lambda_i-\mathcal{Q}_{hi}\lambda_i\|_{0,\Gamma_i}$. Then summing \eqref{eq:staqhi} from $i=1$ to $i=K$ leads to \eqref{eq:staqh}. This finishes the proof.
\end{proof}

For $\lambda_h\in\mathcal{B}_h$, denote its decomposition by
$\lambda_h=\mathring{\lambda}_h+\lambda_h^c$, with
$\mathring{\lambda}_h\in\mathring{\mathcal{B}}_h$ and
$\lambda_h^c\in\mathcal{B}_h^c$, such that $\mathring{\lambda}_h(x)=0$
if $x$ is a corner point, and $\lambda_h^c(x)=0$ if $x$ is a vertex
other than the corners.

\begin{lemma}
For any $\lambda_h\in\mathcal{B}_h$, $\|\mathring{\lambda}_h\|_{0,\Gamma}\lesssim h^{-1/2}\|\mathring{\lambda}_h\|_{-1/2,c,\Gamma}$ and $\|\lambda_h^c\|_{-1/2,c,\Gamma}\lesssim h^{1/2}\|\lambda_h^c\|_{0,\Gamma}$.
\end{lemma}

\begin{proof}
By definition of $H^{1/2}_c(\Gamma)$ and inverse inequality, we have
$$
\begin{array}{rl}
\displaystyle\|\mathring{\lambda}_h\|_{0,\Gamma}^2=(\mathring{\lambda}_h,\mathring{\lambda}_h)_\Gamma&\displaystyle\leqslant \|\mathring{\lambda}_h\|_{-1/2,c,\Gamma}\|\mathring{\lambda}_h\|_{1/2,c,\Gamma}\cequiv \|\mathring{\lambda}_h\|_{-1/2,c,\Gamma}\sum_{i=1}^K\|\mathring{\lambda}_h\chi_i\|_{1/2,\Gamma}\\
&\displaystyle\lesssim \|\mathring{\lambda}_h\|_{-1/2,c,\Gamma}\sum_{i=1}^K(h^{-1/2}\|\mathring{\lambda}_h\chi_i\|_{0,\Gamma})\cequiv h^{-1/2}\|\mathring{\lambda}_h\|_{-1/2,c,\Gamma}\|\mathring{\lambda}_h\|_{0,\Gamma}.
\end{array}
$$
Thus, $\|\mathring{\lambda}_h\|_{0,\Gamma}\lesssim h^{-1/2}\|\mathring{\lambda}_h\|_{-1/2,c,\Gamma}$. To prove the second part, for any $\varphi\in H^{1/2}_c(\Gamma)$, denote $\varphi_h=\mathcal{Q}_h\varphi$ and $(\lambda_h^c,\varphi)_\Gamma=(\lambda_h^c,\varphi_h)_\Gamma+(\lambda_h^c,\varphi-\varphi_h)_\Gamma\lesssim h^{1/2}\|\lambda_h^c\|_{0,\Gamma}\|\varphi_h\|_{1/2,c,\Gamma}+h^{1/2}\|\lambda_h^c\|_{0,\Gamma}\|\varphi\|_{1/2,c,\Gamma}$. Therefore, $\|\lambda_h^c\|_{-1/2,c,\Gamma}=\sup_{\varphi\in H^{1/2}_c(\Gamma)}\frac{(\lambda_h^c,\varphi)}{\|\varphi\|_{1/2,c,\Gamma}}\lesssim h^{1/2}\|\lambda_h^c\|_{0,\Gamma}.$ This finishes the proof.
\end{proof}

\begin{lemma}\label{lem:sdcorner}
It holds for $\lambda_h\in\mathcal{B}_h$ that
$$
\|\lambda_h\|_{-1/2,c,\Gamma}+h^{1/2}\|\lambda_h\|_{0,\Gamma}\cequiv \|\mathring{\lambda}_h\|_{-1/2,c,\Gamma}+h^{1/2}\|\lambda_h^c\|_{0,\Gamma}.
$$
\end{lemma}
\begin{proof}
On one hand,
$$
\begin{array}{rl}
\displaystyle\|\lambda_h\|_{-1/2,c,\Gamma}+h^{1/2}\|\lambda_h\|_{0,\Gamma} & \displaystyle\leqslant  \|\mathring{\lambda}_h\|_{-1/2,c,\Gamma}+h^{1/2}\|\mathring{\lambda}_h\|_{0,\Gamma}+\|\lambda_h^c\|_{-1/2,c,\Gamma}+h^{1/2}\|\lambda_h^c\|_{0,\Gamma}\\
&\displaystyle\lesssim \|\mathring{\lambda}_h\|_{-1/2,c,\Gamma}+h^{1/2}\|\lambda_h^c\|_{0,\Gamma};
\end{array}
$$
on the other hand,
$$
\begin{array}{rl}
\displaystyle\|\mathring{\lambda}_h\|_{-1/2,c,\Gamma}+h^{1/2}\|\lambda_h^c\|_{0,\Gamma} & \displaystyle\leqslant \|\lambda_h\|_{-1/2,c,\Gamma}+\|\lambda_h^c\|_{-1/2,c,\Gamma}+h^{1/2}\|\lambda_h^c\|_{0,\Gamma}\\
& \displaystyle\lesssim \|\lambda_h\|_{-1/2,c,\Gamma}+h^{1/2}\|\lambda_h^c\|_{0,\Gamma}\lesssim \|\lambda_h\|_{-1/2,c,\Gamma}+h^{1/2}\|\lambda_h\|_{0,\Gamma}.
\end{array}
$$
The proof is thus finished.
\end{proof}

\begin{lemma}\label{lem:locanegnorm}
It holds for $\lambda_h\in\mathcal{B}_h$ that
\begin{equation}
\|\lambda_h\|_{-1/2,c,\Gamma}+h^{1/2}\|\lambda_h\|_{0,\Gamma}\lesssim \|\lambda_h\|_{-1/2,\Gamma}\lesssim \|\mathring{\lambda}_h\|_{-1/2,\Gamma}+(1+|\log h|)h^{1/2}\|\lambda_h^c\|_{0,\Gamma}.
\end{equation}
\end{lemma}
\begin{proof}
By inverse inequality, $h^{1/2}\|\lambda_h\|_{0,\Gamma}\lesssim \|\lambda_h\|_{-1/2,\Gamma}$. By definition,
\begin{equation}\label{eq:negcomp}
\|\lambda_h\|_{-1/2,c,\Gamma}=\sup_{\mu\in H^{1/2}_c(\Gamma)}\frac{(\lambda_h,\mu)_\Gamma}{\|\mu\|_{1/2,c,\Gamma}}\lesssim \sup_{\mu\in H^{1/2}_c(\Gamma)}\frac{(\lambda_h,\mu)_\Gamma}{\|\mu\|_{1/2,\Gamma}}\lesssim \sup_{\mu\in H^{1/2}(\Gamma)}\frac{(\lambda_h,\mu)_\Gamma}{\|\mu\|_{1/2,\Gamma}}=\|\lambda_h\|_{-1/2,\Gamma}.
\end{equation}
The left inequality is then proved. For the right one, we only have to note that $\|\lambda_h^c\|_{-1/2,\Gamma}\lesssim h^{1/2}(1+|\log h|)\|\lambda_h^c\|_{0,\Gamma}$, and the inequality follows from the triangle inequality.
\end{proof}

\begin{lemma}\label{lem:equinormint}
It holds for $\lambda_h\in\mathring{\mathcal{B}_h}$ that
\begin{equation}\label{eq:equihalfnorm}
\|\lambda_h\|_{1/2,\Gamma}\lesssim \|\lambda_h\|_{1/2,c,\Gamma} \lesssim(1+|\log h|)\|\lambda_h\|_{1/2,\Gamma}
\end{equation}
and
\begin{equation}\label{eq:equihalfnorm2}
\|\lambda_h\|_{-1/2,\Gamma}\lesssim (1+|\log h|)^{2}\|\lambda_h\|_{-1/2,c,\Gamma}.
\end{equation}
\end{lemma}
\begin{proof}
For any $\lambda_h\in\mathring{\mathcal{B}}_h$, $\|\lambda_h\|_{1/2,\Gamma}\leqslant\sum_{i=1}^K\|\lambda_h\chi_i\|_{1/2,\Gamma}\cequiv\|\lambda_h\|_{1/2,c,\Gamma}$ and $\|\lambda_h\|_{1/2,c,\Gamma}\cequiv\sum_{i=1}^K\|\lambda_h\chi_i\|_{1/2,\Gamma}\lesssim (1+|\log h|)\|\lambda_h\|_{1/2,\Gamma}$. Thus, \eqref{eq:equihalfnorm} is proved.

We now turn to \eqref{eq:equihalfnorm2}. For any $\mu_h\in\mathcal{B}_h$, by Lemma \ref{lem:-12} and \eqref{eq:equihalfnorm}, $\|\mathring{\mu}_h\|_{1/2,c,\Gamma}+\|\mu_h^c\|_{1/2,\Gamma} \lesssim (1+|\log h|)^{2} \|\mu_h\|_{1/2,\Gamma}$.
Therefore,
$$
\begin{array}{rl}
\displaystyle(\lambda_h,\mu_h)_\Gamma &\displaystyle\leqslant \|\lambda_h\|_{-1/2,c,\Gamma}\|\mathring{\mu}_h\|_{1/2,c,\Gamma}+\|\lambda_h\|_{0,\Gamma}\|\mu_h^c\|_{0,\Gamma}\lesssim \|\lambda_h\|_{-1/2,c,\Gamma}\|\mathring{\mu}_h\|_{1/2,c,\Gamma}+h^{1/2}\|\lambda_h\|_{0,\Gamma}\|\mu_h^c\|_{1/2,\Gamma}\\
&\displaystyle\lesssim \|\lambda_h\|_{-1/2,c,\Gamma}(\|\mathring{\mu}_h\|_{1/2,c,\Gamma}+\|\mu_h^c\|_{1/2,\Gamma})\lesssim (1+|\log h|)^2\|\lambda_h\|_{-1/2,c,\Gamma}\|\mu_h\|_{1/2,\Gamma},
\end{array}
$$
and $\|\lambda_h\|_{-1/2,\Gamma}\cequiv\sup_{\mu_h\in\mathcal{B}_h}\frac{(\lambda_h,\mu_h)_\Gamma}{\|\mu_h\|_{1/2,\Gamma}}\lesssim (1+|\log h|)^2\|\lambda_h\|_{-1/2,c,\Gamma}$. This finishes the proof.
\end{proof}

\subsection{Discrete harmonic extension operator and the discrete inverse generalized Poincar\'e--Steklov operator}\label{sec:dishedisgips}

For $\lambda_h\in\mathcal{B}_h$, let $u_{h,\lambda_h}$ satisfy 
\begin{equation}
\left\{\begin{array}{ll}
\displaystyle (\nabla u_{h,\lambda_h},\nabla v_h)=0,&\forall\,v_h\in V_{h0},\\
\displaystyle u_{h,\lambda_h}=\lambda_h,&\mbox{on}\,\partial\Omega.
\end{array}
\right.
\end{equation}
Define the discrete harmonic operator $E_h:\mathcal{B}_h\to V_h$ by $E_h\lambda_h=u_{h,\lambda_h}$. Simultaneously, $(E_h\lambda_h,\Delta_{h,1}v_h)=0$ for $v_h\in V_{h0}$. The stability properties of $E$ and $E_c$ can be inherited by $E_h$. The lemma below can be referred to \cite{GlowinskiPironneau1979}.

\begin{lemma}\label{lem:cghbi} $\|E_h\lambda_h\|_{1,\Omega}\cequiv \|\lambda_h\|_{1/2,\Gamma}$ for $\lambda_h\in\mathcal{B}_h$.
\end{lemma}

Define  $\widetilde{\mathcal{B}_h}:=\Big\{\lambda_h\in\mathcal{B}_h:E_h\lambda_h\in\mathcal{R}(\Delta)\Big\}$. Then, $E_h\widetilde{\mathcal{B}}_h=E_h\mathcal{B}_h\cap\mathcal{R}(\Delta)$, and by Lemma \ref{lem:codim}, $codim(\widetilde{\mathcal{B}_h},\mathcal{B}_h)=codim(E_h\widetilde{B}_h,E_h\mathcal{B}_h)=codim(E_h\mathcal{B}_h\cap\mathcal{R}(\Delta),E_h\mathcal{B}_h)\leqslant m_0$.

\begin{lemma}\label{lem:bihesta}
The following stability results hold:
\begin{enumerate}
\item it holds for $\lambda_h\in\mathcal{B}_h$ that
$$
\|E_h\lambda_h\|_{0,\Omega}\gtrsim \|\lambda_h\|_{-1/2,c}+h^{1/2}\|\lambda_h\|_{0,\Gamma};
$$

\item it holds for $\lambda_h\in\widetilde{\mathcal{B}_h}$ that
$$
\|E_h\lambda_h\|_{0,\Omega}\lesssim\|\lambda_h\|_{-1/2,c}+h^{1/2}\|\lambda_h\|_{0,\Gamma}.
$$
\end{enumerate}
\end{lemma}

\begin{proof}
Given $\lambda_h\in\mathcal{B}_h$, let $v\in H^2(\Omega)\cap H^1_0(\Omega)$, then
\begin{eqnarray*}
|(E_c\lambda_h-E_h\lambda_h,\Delta v)|
&=&|(\lambda_h,\frac{\partial v}{\partial\mathbf{n}})_\Gamma-(E_h\lambda_h,\Delta v)|
=|(\nabla E_h\lambda_h,\nabla v)|\\
&=&\inf_{v_h\in V_{h0}}|(\nabla E_h\lambda_h,\nabla (v-v_h))|
\lesssim h\|E_h\lambda_h\|_{1,\Omega}|v|_{2,\Omega}\\
&\lesssim& h^{1/2}\|\lambda_h\|_{0,\Gamma}\|\Delta v\|_{0,\Omega}.
\end{eqnarray*}

First, taking $v$ such that $\Delta v=E_c\lambda_h$ gives
$$
\|E_h\lambda_h\|_{0,\Omega}
\geqslant c_1\|\lambda_h\|_{-1/2,c}-c_2h^{1/2}\|\lambda_h\|_{0,\Gamma},
$$
and further, $\|E_h\lambda_h\|_{0,\Omega}\geqslant c_3h^{1/2}\|\lambda_h\|_{0,\Gamma}$. It is obtained that
$$
\|E_h\lambda_h\|_{0,\Omega}\geqslant\max(\frac{c_1c_3}{c_2+c_3}\|\lambda_h\|_{-1/2,c},c_3h^{1/2}\|\lambda_h\|_{0,\Gamma}).
$$

For $\lambda_h\in\mathcal{B}_h$ such that $E_h\lambda_h\in\mathcal{R}(\Delta)$,
taking $v$ such that $\Delta v=E_h\lambda_h$ yields
$$
\|E_h\lambda_h\|_{0,\Omega}\lesssim \|\lambda_h\|_{-1/2,c}+h^{1/2}\|\lambda_h\|_{0,\Gamma}.
$$

The lemma is proved.
\end{proof}

Define the discrete inverse generalized Poincar\'e--Steklov operator $S_h:\mathcal{B}_h\to \mathcal{B}_h$ by
\begin{equation}\label{eq:defgh}
(\lambda_h,S_h\gamma_h)_\Gamma=(E_h\lambda_h,E_h\gamma_h),\ \forall\,\lambda_h,\gamma_h\in\mathcal{B}_h.
\end{equation}
The decomposition of the first biharmonic problem can be inherited in the discrete level.
\begin{proposition}
Let $u_h\in V_{h0}$ solve $(\Delta_{h,1} u_h,\Delta_{h,1} v_h)=(f,v_h)$, $\forall\, v_h\in V_{h0}$, then $u_h$ can be obtained by seeking $(\tilde{u}_h,\zeta_h,u_h)\in V_{h0}\times\mathcal{B}_h\times V_{h0}$, such that
\begin{eqnarray}
&(\Delta_{h,2} \tilde{u}_h,\Delta_{h,2} p_h)=(f,p_h),&\forall\, p_h\in V_{h0};\label{dgpdis1}\\
&(\zeta_h,S_h\gamma_h)_\Gamma=-(\Delta_{h,2}\tilde{u}_h,E_h\gamma_h),&\forall\,\gamma_h\in\mathcal{B}_h;\label{dgpdis2}\\
&(\Delta_{h,2} u_h,\Delta_{h,2} v_h)=(\Delta_{h,2}\tilde{u}_h+E_h\zeta_h,\Delta_{h,2} v_h),&\forall\,v_h\in V_{h0}.\label{dgpdis3}
\end{eqnarray}
\end{proposition}

Define $\mathcal{F}_h:\mathcal{B}_h\to \mathcal{B}_h$ such that for  $\gamma_h,\mu_h\in\mathcal{B}_h$,
\begin{equation}\label{eqn:th}
(\mathcal{F}_h\gamma_h,\mu_h)_\Gamma: =(\nabla E_h\gamma_h,\nabla E_h\mu_h)+(E_h\gamma_h,E_h\mu_h),
\end{equation}
and $\mathcal{D}_h:\mathring{\mathcal{B}}_h\to \mathring{\mathcal{B}}_h$ such that for $\gamma_h,\mu_h\in\mathring{\mathcal{B}}_h$
\begin{equation}
(\mathcal{D}_h\gamma_h,\mu_h)_\Gamma:=\sum_{i=1}^K(\nabla E_h(\chi_i\gamma_h),\nabla E_h(\chi_i\mu_h)).
\end{equation}
\begin{lemma}\label{lem:bdycnorm}
$\mathcal{F}_h$ is an isomorphism from $(\mathcal{B}_h,\|\cdot\|_{1/2,\Gamma})$ onto $(\mathcal{B}_h,\|\cdot\|_{-1/2,\Gamma})$, and $\mathcal{D}_h$ is an isomorphism from $(\mathring{\mathcal{B}}_h,\|\cdot\|_{1/2,c,\Gamma})$ onto $(\mathring{\mathcal{B}}_h,\|\cdot\|_{-1/2,c,\Gamma})$.
\end{lemma}
Some equivalent description of $S_h$ can be established by means of $\mathcal{D}_h$ and $\mathcal{F}_h$.
\begin{theorem}\label{thm:halfchi}
For any $\lambda_h\in\mathcal{B}_h$, it holds that $$
(\lambda_h,S_h\lambda_h)_\Gamma\gtrsim (\mathcal{D}_h^{-1}\mathring{\lambda}_h,\mathring{\lambda}_h)_\Gamma+h(\lambda_h^c,\lambda_h^c)_\Gamma,
$$
and for any $\lambda_h\in\widetilde{\mathcal{B}_h}$,
$$
(\lambda_h,S_h\lambda_h)_\Gamma\lesssim (\mathcal{D}_h^{-1}\mathring{\lambda}_h,\mathring{\lambda}_h)_\Gamma+h(\lambda_h^c,\lambda_h^c)_\Gamma.
$$
\end{theorem}
\begin{proof}
The theorem follows from Lemmas \ref{lem:bihesta}, \ref{lem:sdcorner}, and \ref{lem:bdycnorm} directly.
\end{proof}

\begin{theorem}\label{thm:sdbihe}
For any $\lambda_h\in\mathcal{B}_h$, it holds that $$
(\lambda_h,S_h\lambda_h)_\Gamma\gtrsim (1+|\log h|)^{-4}\Big((\mathcal{F}_h^{-1}\mathring{\lambda}_h,\mathring{\lambda}_h)_\Gamma+h(\lambda_h^c,\lambda_h^c)_\Gamma\Big),
$$
and for any $\lambda_h\in\widetilde{\mathcal{B}_h}$,
$$
(\lambda_h,S_h\lambda_h)_\Gamma\lesssim (\mathcal{F}_h^{-1}\mathring{\lambda}_h,\mathring{\lambda}_h)_\Gamma+h(\lambda_h^c,\lambda_h^c)_\Gamma.
$$
\end{theorem}
\begin{proof}
The theorem follows from Lemmas \ref{lem:equinormint}, \ref{lem:sdcorner}, and \ref{lem:bdycnorm} directly.
\end{proof}

\begin{theorem}\label{thm:sdbihe14}
For any $\lambda_h\in\mathcal{B}_h$, it holds that
$$
(\lambda_h,S_h\lambda_h)_\Gamma\gtrsim (1+|\log h|)^{-4}(\mathcal{F}_h^{-1}\lambda_h,\lambda_h)_\Gamma,
$$
and for any $\lambda_h\in\widetilde{\mathcal{B}_h}$, it holds that
$$
(\lambda_h,S_h\lambda_h)_\Gamma\lesssim (\mathcal{F}_h^{-1}\lambda_h,\lambda_h)_\Gamma.
$$
\end{theorem}
\begin{proof}
The theorem follows from Theorem \ref{thm:sdbihe} and Lemma \ref{lem:locanegnorm}.
\end{proof}
\subsection{A modified stable decomposition of $M_{h,1}$}\label{sec:app:suboptimal}

Define $\widehat{V}_{h0}:=\Big\{p_h\in V_{h0}:(\Delta_{h,1}p_h)|_\Gamma\in\widetilde{\mathcal{B}_h}\Big\}$ and $\widetilde{M}_{h,1}:=\Big\{w_h\in M_{h,1}:I_hw_h\in\widehat{V}_{h0}\Big\}$. Then, by Lemma \ref{lem:codim}, $codim(\widehat{V}_{h0},V_{h0})\leqslant codim(\widetilde{\mathcal{B}}_h,\mathcal{B}_h)\leqslant m_0$ and $codim(\widetilde{M}_{h,1},M_{h,1})=codim(\widehat{V}_{h0}, I_hM_{h,1})\leqslant codim(\widehat{V}_{h0}, V_{h0})\leqslant m_0$.

\begin{lemma}\label{lem:cprdeltas}
  Let $v_h\in V_{h0}$. Let $\Delta_{h,1}$ and $\Delta_{h,2}$ be
  Laplacian operators defined on $V_{h0}$. It holds that
\begin{enumerate}
\item $\|\Delta_{h,2} v_h\|_{0,\Omega}\leqslant \|\Delta_{h,1} v_h\|_{0,\Omega}$;
\item $\|\Delta_{h,1}v_h\|_{0,\Omega}\lesssim h^{-1/2}\|\Delta_{h,2} v_h\|_{0,\Omega}$ for  $v_h\in \widehat{V}_{h0}$.
\end{enumerate}
\end{lemma}
\begin{proof}
Let $q_h\in V_{h0}$, then by the definitions, $(\Delta_{h,1} v_h,q_h)=(\Delta_{h,2}v_h,q_h)$; that is, $\Delta_{h,2}v_h\in V_{h0}$ is the $L^2$ projection of $\Delta_{h,1} v_h$ into $V_{h0}$. The first result follows then.

Now let $\mathcal{H}_h$ be the $H^1-$projection operator from $V_h$ to $V_{h0}$, namely $(\nabla p_h,\nabla q_h)=(\nabla p_h,\nabla \mathcal{H}_hq_h)$, for $\forall\, p_h\in V_{h0}$ and $\forall\,q_h\in V_h$. It follows that $q_h-\mathcal{H}q_h=E_hq_h|_\Gamma$ and $(\Delta_{h,1}p_h,q_h)=(\Delta_{h,2}p_h,\mathcal{H}_hq_h)$. If $q_h|_\Gamma\in \widetilde{\mathcal{B}_h}$, then $\|\mathcal{H}_hq_h-q_h\|_{0,\Omega}\lesssim \|q_h|_\Gamma\|_{-1/2,c,\Gamma}+h^{1/2}\|q_h|_\Gamma\|_{0,\Gamma}\lesssim \|q_h|_\Gamma\|_{0,\Gamma}\lesssim h^{-1/2}\|q_h\|_{0,\Omega},$ and thus $\|\mathcal{H}_hq_h\|_{0,\Omega}\leqslant \|q_h\|_{0,\Omega}+\|\mathcal{H}_hq_h-q_h\|_{0,\Omega}\lesssim h^{-1/2}\|q_h\|_{0,\Omega}$. Therefore for $w_h\in\widehat{V}_{h0}$, $(\Delta_{h,1}w_h,\Delta_{h,1}v_h)=(\Delta_{h,2}w_h,\mathcal{H}_h (\Delta_{h,1}v_h))\lesssim \|\mathcal{H}_h (\Delta_{h,1}v_h)\|_{0,\Omega}\|\Delta_{h,2} w_h\|_{0,\Omega}\lesssim h^{-1/2}\|\Delta_{1,h}v_h\|_{0,\Omega} \|\Delta_{h,2}w_h\|_{0,\Omega}$, and $\|\Delta_{h,1}w_h\|_{0,\Omega}\displaystyle=\sup_{v_h\in \widehat{V}_{h0}}\frac{(\Delta_{h,1}w_h,\Delta_{h,1}v_h)}{\|\Delta_{h,1}v_h\|_{0,\Omega}}\lesssim h^{-1/2}\|\Delta_{h,2}w_h\|_{0,\Omega}$.
The lemma is proved.
\end{proof}

\begin{lemma}\label{lem:mdsd}
The modified stable decomposition holds:
\begin{enumerate}
\item For $w_h\in M_{h,1}$, it holds that
\begin{equation}\label{sdopemd1}
|w_h|_{2,h}^2\gtrsim \sum_{T\in\mathcal{T}_h} h_T^{-4}\|w_h-\Pi_{h,1}I_hw_h\|_{0,T}^2+(\Delta_{h,2}I_hw_h,\Delta_{h,2}I_hw_h).
\end{equation}
\item For $w_h\in\widetilde{M}_{h,1}$, it holds that
\begin{equation}\label{sdopemd2}
|w_h|_{2,h}^2\lesssim \sum_{T\in\mathcal{T}_h} h_T^{-4}\|w_h-\Pi_{h,1}I_hw_h\|_{0,T}^2+h^{-1}(\Delta_{h,2}I_hw_h,\Delta_{h,2}I_hw_h).
\end{equation}
\end{enumerate}
\end{lemma}
\begin{proof}
Combining Theorem \ref{thm:sdcfs} and Lemma \ref{lem:cprdeltas} leads to the lemma directly.
\end{proof}

\section{Optimal solvers for fourth-order finite element problems}\label{sec:opsol}
In this section, we will present several effective solvers for both the first and the second biharmonic problems. We will provide detailed theoretical analysis for these solvers as well as numerical experiments that support our theories.

In the presentation below, $R_h:M_{h,k}'\to M_{h,k}$ represents any symmetric smoother (such as the Jacobi and the symmetric Gauss--Seidel smoother) for the discrete biharmonic operators, and $\Pi_{h,k}$ follows from the definition as in Section \ref{sec:les:fesf}.
\subsection{Second biharmonic finite element problem $A_{h,2}u_h=f_h$ on $M_{h,2}$}\label{sec:ssp}

For the second biharmonic finite element problem, the discrete Laplacian operator of second kind provides an optimal preconditioner. By Theorem \ref{thm:sdcfs}, Lemma \ref{lem:minmax}, and the FASP theory, we obtain the following theorem.

\begin{theorem}\label{thm:opssp}
Define
\begin{equation}
B_{h,2}:=R_h+\Pi_{h,2} \Delta_{h,2}^{-2} (\Pi_{h,2})^*,
\end{equation}
then the $m_0$-th effective condition number $\kappa^{\rm eff}_{m_0}(B_{h,2}A_{h,2})\lesssim 1$.
\end{theorem}

Let $N_p$ be the number of interior points of the triangulation, and let $N_{h,2}$ be the number of the degree of freedoms of $A_{h,2}u_h=f_h$, then $N_{h,2}\cequiv N_p$.

\begin{lemma}\label{lem:GrasedyckXu2011}
\cite{GrasedyckXu2011,Xu2010}
The linear element problem for the Poisson equation can be solved in the complexity of $\mathcal{O}(N_p\log N_p)$.
\end{lemma}
\begin{theorem}
When PCG is applied on $B_{h,2}A_{h,2}u_h=B_{h,2}f_h$, the total complexity is $\mathcal{O}(N_{h,2}\log N_{h,2})$.
\end{theorem}
\begin{proof}
The theorem follows from \eqref{eq:effcondest}, Theorem \ref{thm:opssp} and Lemma \ref{lem:GrasedyckXu2011}.
\end{proof}
\subsection{First biharmonic finite element problem $A_{h,1}u_h=f_h$ on $M_{h,1}$}\label{sec:cpp}
\subsubsection{Preconditioning effect of the discrete Laplacian operator of the second kind}
The discrete Laplacian operator of second kind induces a simple preconditioner for first biharmonic finite element problem.
\begin{theorem}
Define
\begin{equation}\label{FASPR'}
B_{h,1}':=R_h+\Pi_{h,1}\Delta_{h,2}^{-2}(\Pi_{h,1})^*.
\end{equation}
then the $m_0$-th effective condition number $\kappa^{\rm eff}_{m_0}(B_{h,1}'A_{h,1})\lesssim h^{-1}$.
\end{theorem}
\begin{proof}
The theorem follows from Lemma \ref{lem:minmax} and Lemma \ref{lem:mdsd}.
\end{proof}
The computational cost is dominantly contained in Poisson solvers. Let $N_{h,1}$ be the number of the degree of freedoms of $A_{h,1}u_h=f_h$, then $N_{h,1}\cequiv N_{h,2}\cequiv h^{-2}$. 
\begin{theorem}
The complexity of PCG applied on $B_{h,1}'A_{h,1}u_h=B_{h,1}'f_h$ is $\mathcal{O}(N_{h,1}^{1.25}\log N_{h,1})$.
\end{theorem}

\subsubsection{An optimal preconditioner}
\begin{theorem}\label{thm:opcpp}
Define
\begin{equation}\label{FASPreconditioner}
B_{h,1}=R_h+\Pi_{h,1}(\Delta_{h,1}^*\Delta_{h,1})^{-1}(\Pi_{h,1})^*.
\end{equation}
The condition number of $B_{h,1}A_{h,1}$ is bounded uniformly.
\end{theorem}
\begin{proof}
The theorem follows from the FASP theory and the stable decomposition \eqref{sdopefs}.
\end{proof}

The main work of the preconditioner $B_{h,1}$ is in the inversion of $\Delta_{1,h}^*\Delta_{1,h}$: namely, given $f_h\in V_{h0}$, we find $u_h\in V_{h0}$ such that
\begin{equation}\label{femwmodelcppmix}
(\Delta_{h,1}u_h,\Delta_{h,1}v_h)=(f_h,v_h),\quad\forall\,v_h\in V_{h0}.
\end{equation}
\paragraph{\textbf{Strategy of solving \eqref{femwmodelcppmix}}}
Let $u_h$ be the solution of \eqref{femwmodelcppmix}, and $\overline{w}_h:=\Delta_{h,1}u_h$, then $(u_h,\overline{w}_h)$ is the unique solution of the equation
\begin{equation}\label{femfirmix}
\left\{
\begin{array}{ll}
(\nabla u_h,\nabla p_h)=(\overline{w}_h,p_h)&\forall\,p_h\in V_h\\
(\nabla \overline{w}_h,\nabla v_h)=(f,v_h)&\forall\,v_h\in V_{h0}.
\end{array}
\right.
\end{equation}
Namely, \eqref{femwmodelcppmix} is equivalent to \eqref{femfirmix}, whereas the latter is the discretization of the mixed formulation of the first biharmonic problem given in \cite{CiarletRaviart1974}, which finds $(u,\overline{w})\in H^1_0(\Omega)\times H^1(\Omega)$, such that
\begin{equation}\label{mixcont}
\left\{
\begin{array}{ll}
(\nabla u,\nabla v)=(\overline{w},v)&\forall\,v\in H^1(\Omega),\\
(\nabla \overline{w},\nabla p)=(f,p)&\forall\,p\in H^1_0(\Omega).
\end{array}
\right.
\end{equation}

The coupled system \eqref{mixcont} can be decoupled as demonstrated in the following lemma.

\begin{lemma}\label{lem:deccrmixcont}
Let $(u,\overline{w})$ be the solution of \eqref{mixcont}, then they can be obtained by seeking $(u,w,\lambda)\in H^1_0(\Omega)\times H^1_0(\Omega)\times H^{1/2}(\Gamma)$ such that
\begin{eqnarray}
&(\nabla w,\nabla p)=(f,p),&\forall\, p\in H^1_0(\Omega);\label{cdis1}\\
&(E\lambda,E\gamma)=-(w,E\gamma),&\forall\,\gamma\in H^{1/2}(\Gamma);\label{cdis2}\\
&(\nabla u,\nabla v)=(w+E\lambda,v),&\forall\,v\in H^1_0(\Omega).\label{cdis3}
\end{eqnarray}
And, $(u,\overline{w})=(u,w+E\lambda)$.
\end{lemma}
\begin{proof} 
Let $\lambda=\overline{w}|_\Gamma$, then $\overline{w}=w+E\lambda$, with $w\in H^1_0(\Omega)$ uniquely determined by \eqref{cdis1}. Since $(\overline{w},E\gamma)=-(\nabla u,\nabla E\gamma)=0$ for $\gamma\in H^{1/2}(\Gamma)$, $(E\lambda,E\gamma)=-(w,E\gamma)$. Thus, $\lambda$ can be solved from \eqref{cdis2}. Further, $u$ can be solved from \eqref{cdis3}. The proof is finished.
\end{proof}

Analogously, the decoupling of \eqref{femwmodelcppmix} can be carried out by means of $E_h$. Similar to Lemma \ref{lem:deccrmixcont}, we can prove the lemma below.
\begin{lemma}\label{lem:deccrmix}
Let $u_h$ be the solution of \eqref{femwmodelcppmix}, then it can be obtained by seeking $(u_h,w_h,\lambda_h)\in V_{h0}\times V_{h0}\times \mathcal{B}_h$ such that
\begin{eqnarray}
&(\nabla w_h,\nabla p_h)=(f,p_h),&\forall\, p_h\in V_{h0};\label{dis1}\\
&(E_h\lambda_h,E_h\gamma_h)=-(w_h,E_h\gamma_h),&\forall\,\gamma_h\in\mathcal{B}_h;\label{dis2}\\
&(\nabla u_h,\nabla v_h)=(w_h+E_h\lambda_h,v_h),&\forall\,v_h\in V_{h0}.\label{dis3}
\end{eqnarray}
\end{lemma}

According to Lemma \ref{lem:GrasedyckXu2011}, \eqref{dis1} and \eqref{dis3} can both be solved optimally. The difficulty lies in solving \eqref{eq:btransgh} in optimal complexity with respect to the size of \eqref{femwmodelcppmix}. The problem \eqref{dis2} is equivalent to finding $\lambda_h\in\mathcal{B}_h$, such that
\begin{equation}\label{eq:btransgh}
(\lambda_h,S_h\gamma_h)_\Gamma=-(w_h,E_h\gamma_h),\quad\forall\,\gamma_h\in\mathcal{B}_h.
\end{equation}

\paragraph{\textbf{Optimal and nearly optimal preconditioners for $S_h$}}

Based on the equivalent description of $S_h$ in Section \ref{sec:dishedisgips}, we present three preconditioners for $S_h$. They are
\begin{description}
\item[$T_{h,1}$]$=I\mathcal{D}_hI^*+h^{-1}JJ^*$;
\item[$T_{h,2}$]$=I \mathcal{F}_h I^*+h^{-1}J J^*$;
\item[$T_{h,3}$]$=\mathcal{F}_h$.
\end{description}
Here, $I$ is the inclusion operator from $\mathring{\mathcal{B}}_h$ to $\mathcal{B}_h$, $J$ is the inclusion operator from $\mathcal{B}_h^c$ to $\mathcal{B}_h$, and $I^*$ and $J^*$ are the adjoint operators of $I$ and $J$, respectively.

\begin{theorem}\label{thm:opbhesqr}
It holds for $j=1,2,3$ that $\kappa_{m_0}^{\rm eff}(T_{h,j}S_h)\lesssim 1+|\log h|^{\beta_j}$, with $\beta_1=0$, $\beta_2=4$ and $\beta_3=4$.
\end{theorem}
\begin{proof}
The theorem follows directly from the FASP theory, Lemma \ref{lem:minmax}, and Theorems \ref{thm:halfchi}, \ref{thm:sdbihe}, and \ref{thm:sdbihe14}, respectively.
\end{proof}

Let $N_S$ be the number of the degree of freedoms of \eqref{dis2}. Then $N_{h,1}\cequiv N_S^2$, and $N_{h,1}$ is equivalently the number of the degree of freedom of \eqref{femwmodelcppmix}. The theorem below follows from Theorem \ref{thm:opbhesqr}, Lemma \ref{lem:GrasedyckXu2011}, and \eqref{eq:effcondest}.
\begin{theorem}\label{lem:occr}
\begin{enumerate}
\item The equation \eqref{dis2} can be solved by PCG with $T_{h,1}$(or $T_{h,2}$, $T_{h,3}$) as the preconditioner in the complexity of $\mathcal{O}(N_{h,1}\log N_{h,1})$, or $\mathcal{O}(N_{h,1}\log^{3} N_{h,1})$, $\mathcal{O}(N_{h,1}\log^3 N_{h,1})$, respectively.
\item The problem \eqref{femwmodelcppmix} can be solved in the complexity of $\mathcal{O}(N_{h,1}\log N_{h,1})$.
\end{enumerate}
\end{theorem}

Evaluating the action of the operator $(\Delta_{h,1}^*\Delta_{h,1})^{-1}$ requires solving the mixed finite element problem \eqref{femwmodelcppmix} once. Therefore the preconditioner $B_{h,1}$ can be carried out in the complexity of $\mathcal{O}(N_{h,1}\log N_{h,1})$. The theorem below is then derived.
\begin{theorem}\label{thm:occpp}
When PCG is used to solve the problem $B_{h,1}A_{h,1}u_h=B_{h,1}f_h$, the total complexity for convergence is $\mathcal{O}(N_{h,1}\log N_{h,1})$.
\end{theorem}

\section{Numerical examples}\label{sec:numexp}
In this section, we present several numerical examples to illustrate the
effects of the preconditioners given in last section. We will compute
and record the extremal eigenvalues of the preconditioned operators. We will also test the performance of PCG method with given
preconditioners on some modal problems. We run the various  PCG computations with the starting guess $\mathbf{0}$, and with a stop criteria whereby the relative residual ($\|\mbox{residual}\|_{l^2}/\|\mbox{rhs}\|_{l^2}$) is smaller than
$10^{-8}$.

\begin{figure}[htbp]
\begin{minipage}[t]{0.48\linewidth}
\centering
\includegraphics[width=\textwidth]{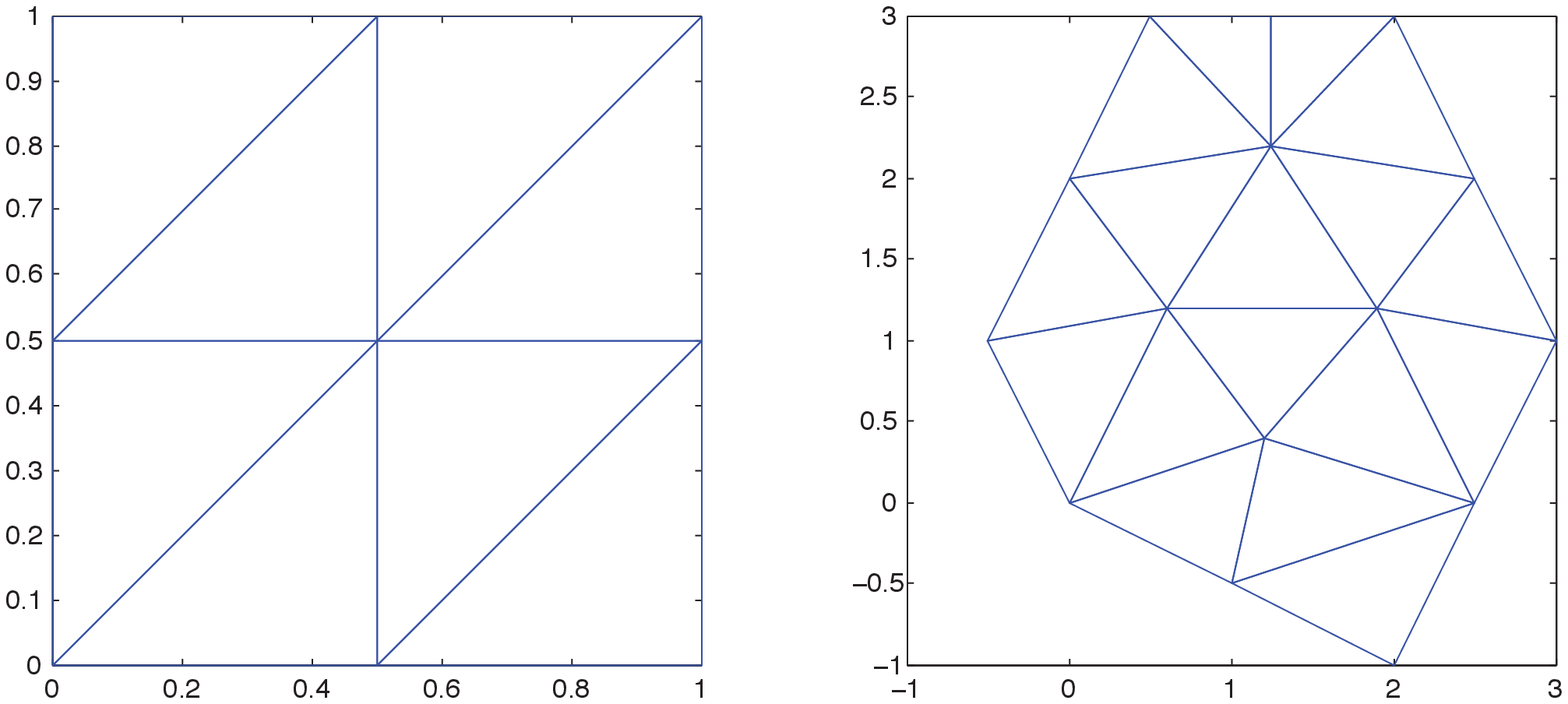}
\end{minipage}
\hfill
\begin{minipage}[t]{0.48\linewidth}
\centering
\includegraphics[width=\textwidth]{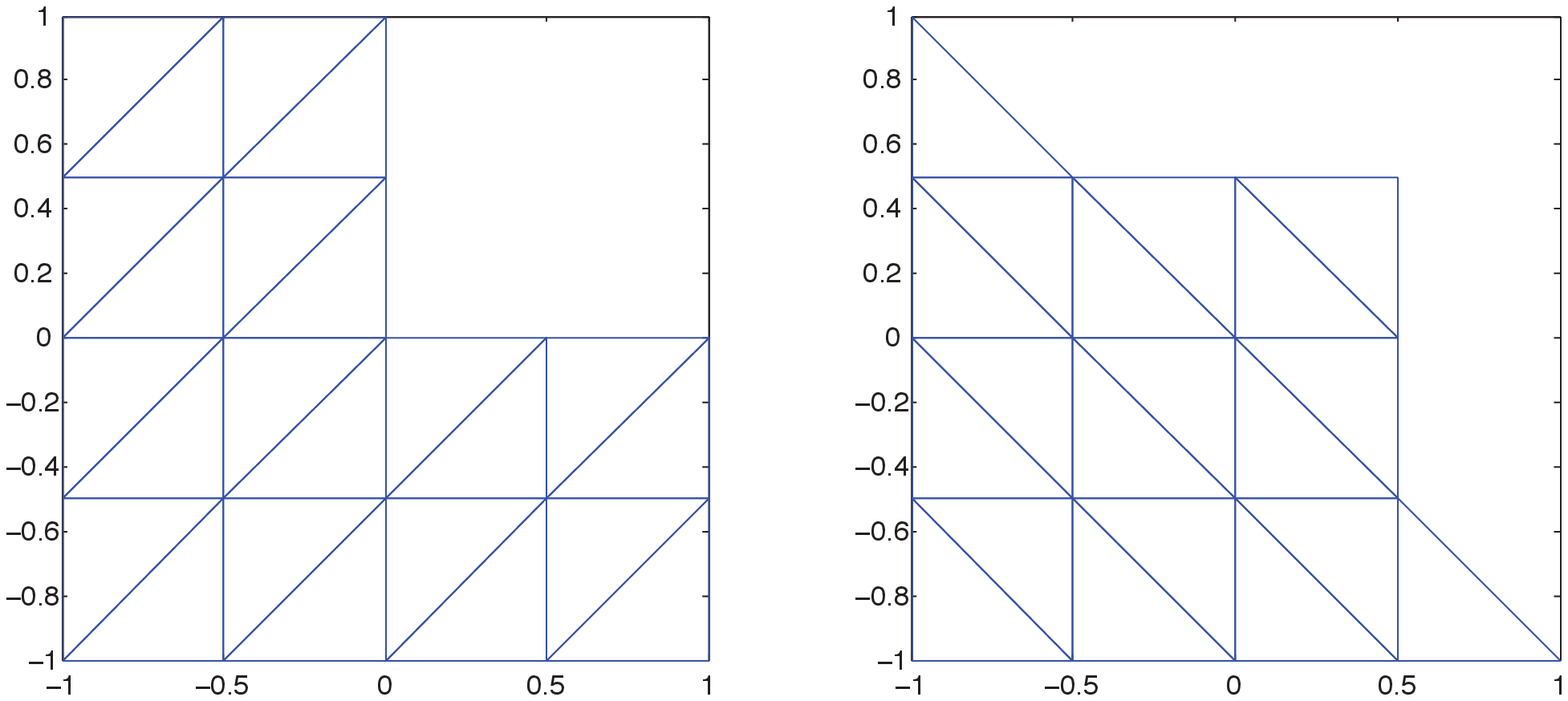}
\end{minipage}
\caption{The computational domains with initial mesh.}\label{fig:initial mesh}
\end{figure}

We test the preconditioner on both convex and nonconvex domains, as shown in Figure \ref{fig:initial mesh}. We divide each computational domain by successively refined quasi-uniform meshes, and we carry out numerical experiments on the multiple meshes to test the performance of each preconditioner. In the tables below, we use $\lambda$ for an eigenvalue, $\kappa$ for a condition number, and DOF for the number of the degree of freedom.

\subsection{Preconditioning effect of $B_{h,1}$ for first biharmonic finite element problem}

We take the frequently used Morley element (\cite{Morley1968}) and the modified Zienkiewicz element (\cite{WangShiXu2007}) as examples to demonstrate the preconditioning effect of the preconditioner $B_{h,1}$ with respect to $A_{h,1}$. We use the symmetric Gauss--Seidel method with three iterations as a smoother.

As shown in Tables \ref{tab:expI}, \ref{tab:expII}, \ref{tab:expIII}, and \ref{tab:expIV}, we observe that for both the Morley element and the modified Zienkiewicz element and on both convex and nonconvex computational domains, the condition number of $B_{h,1}A_{h,1}$ is uniformly bounded as the triangulation is refined.

\begin{table}[h!!]
\centering\small
\begin{tabular}{|c||c|c|c|c||c|c|c|c|}
\hline
\multirow{2}{*}{nodes}& \multicolumn{4}{c||}{Morley Element}&\multicolumn{4}{|c|}{Modified Zienkiewicz Element}\\
\cline{2-9} &DOF&$\lambda_1$&$\lambda_{N_{h,1}}$&$\kappa$&DOF&$\lambda_{1}$&$\lambda_{N_{h,1}}$&$\kappa$
\\
\hline 961&3969&0.70&1.61&2.28&2883&0.18&1.82&10.08\\
\hline 3969&16129&0.70&1.61&2.30&11907&0.18&1.82&10.06\\
\hline 16129&65025&0.70&1.61&2.29&48387&0.18&1.82&10.12\\
\hline
\end{tabular}
\vskip 0.2cm \caption{Eigenvalues of $B_{h,1}A_{h,1}$: unit square.} \label{tab:expI}
\end{table}

\begin{table}[h!!]
\centering\small
\begin{tabular}{|c||c|c|c|c||c|c|c|c|}
\hline
\multirow{2}{*}{nodes}& \multicolumn{4}{c||}{Morley Element}&\multicolumn{4}{|c|}{Modified Zienkiewicz Element}\\
\cline{2-9} &DOF&$\lambda_{1}$&$\lambda_{N_{h,1}}$&$\kappa$&DOF&$\lambda_{1}$&$\lambda_{N_{h,1}}$&$\kappa$\\
\hline 2089&8529&0.70&1.61&2.30&6267&0.16&1.74&10.68\\
\hline 8529&34465&0.70&1.60&2.29&25587&0.16&1.74&10.83\\
\hline 34465&138561&0.70&1.61&2.29&103395&0.16&1.74&10.90\\
\hline
\end{tabular}
\vskip 0.2cm \caption{Eigenvalues of $B_{h,1}A_{h,1}$: hexagon.} \label{tab:expII}
\end{table}

\begin{table}[h!!]
\centering\small
\begin{tabular}{|c||c|c|c|c||c|c|c|c|}
\hline
\multirow{2}{*}{nodes}& \multicolumn{4}{c||}{Morley Element}&\multicolumn{4}{|c|}{Modified Zienkiewicz Element}\\
\cline{2-9} &DOF&$\lambda_{1}$&$\lambda_{N_{h,1}}$&$\kappa$&DOF&$\lambda_{1}$&$\lambda_{N_{h,1}}$&$\kappa$\\
\hline 3201&12033&0.57&1.61&2.83&8835&0.16&1.68&10.55\\
\hline 12545&48641&0.56&1.61&2.86&36099&0.16&1.68&10.53\\
\hline 49665&195585&0.56&1.61&2.88&145923&0.16&1.69&10.55\\
\hline
\end{tabular}
\vskip 0.2cm \caption{Eigenvalues of $B_{h,1}A_{h,1}$: ``L-''shape domain.} \label{tab:expIII}
\end{table}

\begin{table}[h!!]
\centering\small
\begin{tabular}{|c||c|c|c|c||c|c|c|c|}
\hline
\multirow{2}{*}{nodes}& \multicolumn{4}{c||}{Morley Element}&\multicolumn{4}{|c|}{Modified Zienkiewicz Element}\\
\cline{2-9} &DOF&$\lambda_{1}$&$\lambda_{N_{h,1}}$&$\kappa$&DOF&$\lambda_{1}$&$\lambda_{N_{h,1}}$&$\kappa$\\
\hline 2673&10017&0.70&1.61&2.29&7347&0.16&1.65&10.31\\
\hline 10465&40513&0.70&1.61&2.30&30051&0.16&1.65&10.34\\
\hline 41409&162945&0.70&1.61&2.30&1539&0.16&1.66&10.35\\
\hline
\end{tabular}
\vskip 0.2cm \caption{Eigenvalues of $B_{h,1}A_{h,1}$: trident  domain.} \label{tab:expIV}
\end{table}

\subsection{Preconditioning effect of $B_{h,2}$ for the second biharmonic finite element problem.}

We take the Morley element (\cite{Morley1968},\cite{VeigaNiiranenStenberg2010}) to demonstrate the preconditioning effect of the preconditioner $B_{h,2}$ with respect to $A_{h,2}$.  Numerical verifications (\cite{VeigaNiiranenStenberg2010}) showed that the Morley element fits the second biharmonic problem.

As shown in Tables \ref{tab:ssp:pcgstep}, \ref{tab:ssp:convex}, and \ref{fig:ssp:eev}, we observe that for both convex and nonconvex cases the number of PCG iterations remains nearly constant as the number of DOF grows. Moreover, when the computational domain is convex, the eigenvalues of the preconditioned operator are bounded uniformly. However, when the computational domain is not convex, the eigenvalues of the preconditioned operator are bounded uniformly from below. In addition, when the domain is not convex, the eigenvalues are bounded uniformly from above with no more than $m_0$ exceptions.

\begin{table}[htbp]
\centering\small
\begin{tabular}{|c||c|c||c|c||c|c||c|c|}
\hline
Refine & \multicolumn{2}{c||}{Four-Square}&\multicolumn{2}{|c||}{Hexagon}&\multicolumn{2}{|c||}{``L-''shape domain}&\multicolumn{2}{|c|}{Trident}\\
\cline{2-9}  Times&DOF&  step&DOF&  step&DOF&  step&DOF& step\\
\hline 5&16385&14&34817&11&49153&17&40961&12\\
\hline 6&65537&14&139265&11&196609&17&163841&12\\
\hline 7&262145&15&557057&10&786433&19&655361&12\\
\hline
\end{tabular}
\vskip 0.2cm \caption{PCG steps needed for solving $B_{h,2}A_{h,2}$.} \label{tab:ssp:pcgstep}
\end{table}

\begin{table}[htbp]
\centering\small
\begin{tabular}{|c||c|c|c|c||c|c|c|c|}
\hline
Refine & \multicolumn{4}{c||}{Four-Square}&\multicolumn{4}{|c|}{Hexagon}\\
\cline{2-9}  Times&DOF&$\lambda_{1}$&$\lambda_{N_{h,2}}$&$\kappa$&DOF&$\lambda_{1}$&$\lambda_{N_{h,2}}$&$\kappa$\\
\hline 5&16385&0.71&1.60&2.25&34817&0.75&1.57&2.09\\
\hline 6&65537&0.71&1.60&2.25&139265&0.78&1.53&1.96\\
\hline 7&262145&0.71&1.60&2.25&557057&0.81&1.51&1.86\\
\hline
\end{tabular}
\vskip 0.2cm \caption{Eigenvalues of $B_{h,2}A_{h,2}$: convex domain.} \label{tab:ssp:convex}
\end{table}

\begin{table}[htbp]
\centering\small
\begin{tabular}{|c|c|c|c|c||c|c|c|c|c|c|}
\hline
 \multicolumn{5}{|c||}{``L-"shape Domain($m_0=1$)}&\multicolumn{6}{|c|}{Trident Domain($m_0=2$)}\\
\hline DOF&$\lambda_{1}$&$\lambda_{N_S-1}$&$\lambda_{N_S}$&$\kappa^{\rm eff}_{m_0}$&DOF&$\lambda_{min}$&$\lambda_{N_S-2}$&$\lambda_{N_S-1}$&$\lambda_{N_S}$&$\kappa^{\rm eff}_{m_0}$\\
\hline 193&0.73&1.61&1.87&2.22&161&0.69&1.61&1.64&1.68&2.34\\
\hline 769&0.69&1.61&2.33&2.33&641&0.70&1.61&1.67&1.72&2.29\\
\hline 3073&0.71&1.62&3.17&2.28&2561&0.71&1.62&1.80&1.88&2.28\\
\hline
\end{tabular}
\vskip 0.2cm \caption{Eigenvalues of $B_{h,2}A_{h,2}$: nonconvex domain.} \label{fig:ssp:eev}
\end{table}

\subsection{Preconditioning effects of the preconditioners for $S_h$}

\paragraph{\textbf{The optimal preconditioner $T_{h,1}$}} We record the performance of $T_{h,1}$ for $S_h$ in Table \ref{tab:expt1convex} and Table  \ref{tab:expt1nonconvex}. We observe that, the eigenvalues of $T_{h,1}S_h$ are uniformly bounded from below. In addition, they are bounded from above with no more than $m_0$ exceptions.

\begin{table}[htbp]
\centering\small
\begin{tabular}{|c|c|c|c||c|c|c|c|}
\hline
 \multicolumn{4}{|c||}{Unit Square}&\multicolumn{4}{|c|}{Convex Hexagon}\\
\hline DOF&$\lambda_{1}$&$\lambda_{N_S}$&$\kappa$&DOF&$\lambda_{1}$&$\lambda_{N_S}$&$\kappa$\\
\hline 256&0.20&4.48&21.96&352&0.21&4.87&23.40\\
\hline 512&0.20&4.49&22.52&704&0.20&4.89&24.13\\
\hline 1024&0.20&4.49&22.86&1408&0.20&4.90&24.62\\
\hline
\end{tabular}
\vskip 0.2cm \caption{Eigenvalues of $T_{h,1}S_h$: convex domain.} \label{tab:expt1convex}
\end{table}

\begin{table}[htbp]
\centering\small
\begin{tabular}{|c|c|c|c|c||c|c|c|c|c|c|}
\hline
 \multicolumn{5}{|c||}{``L-"shape Domain($m_0=1$)}&\multicolumn{6}{|c|}{Trident Domain($m_0=2$)}\\
\hline DOF&$\lambda_{1}$&$\lambda_{N_S-1}$&$\lambda_{N_S}$&$\kappa^{\rm eff}_{m_0}$&DOF&$\lambda_{min}$&$\lambda_{N_S-2}$&$\lambda_{N_S-1}$&$\lambda_{N_S}$&$\kappa^{\rm eff}_{m_0}$\\
\hline 512&0.21&4.50&35.37&21.89&448&0.19&4.50 &14.70&17.08&24.07\\
\hline 1024&0.20&4.50&56.70&22.39&896&0.18 &4.50&20.24 &23.39 &24.63\\
\hline 2048&0.20&4.50&90.58& 22.76&1792&0.18&4.50 &27.59& 31.75&25.01\\
\hline
\end{tabular}
\vskip 0.2cm \caption{Eigenvalues of $T_{h,1}S_h$: nonconvex domain.} \label{tab:expt1nonconvex}
\end{table}

\paragraph{\textbf{The nearly optimal preconditioners $T_{h,2}$ and $T_{h,3}$}}

We also carry out the same numerical experiments for the performance of $T_{h,2}$ and $T_{h,3}$ as preconditioners for $S_h$, recording the results in Tables \ref{tab:expV},\ref{tab:expVI},\ref{tab:expVII} and \ref{tab:expVIII}. We observe that the eigenvalues of $T_{h,2}S_h$ and $T_{h,3}S_h$ are bounded from above, with no more than $m_0$ exceptions. In addition, they are bounded from below, with slight decrease.

\begin{figure}[htbp]
\includegraphics[width=\textwidth]{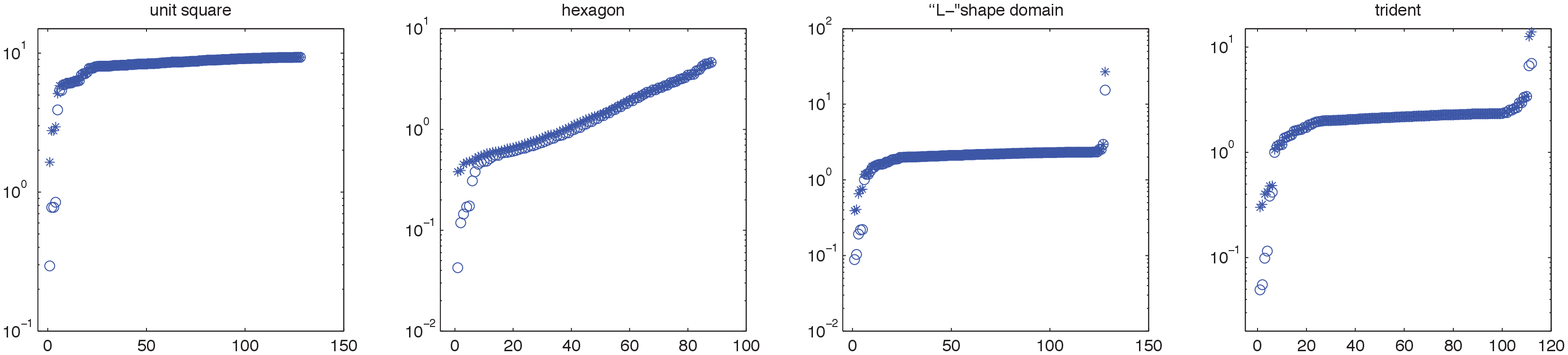} 
\vskip 0cm\caption{The distributions of eigenvalues of $T_{h,2}S_h$ and $T_{h,3}S_h$. In each subfigure, we use $\divideontimes$ for the eigenvalue of $T_{h,2}S_h$, and $\circ$ for that of $T_{h,3}S_h$.}\label{fig:compareQB}
\end{figure}

To illustrate the difference and make a comparison, we compute the eigenvalues of $T_{h,2}S_h$ and $T_{h,3}S_h$ on all the computational domains and plot the distribution of them in Figure \ref{fig:compareQB}. We observe in Figure \ref{fig:compareQB} that, $T_{h,2}$ performs better than $T_{h,3}$ at capturing the extremely low-frequency part, and is as good as $T_{h,3}$ at capturing the high-frequency part.

\begin{table}[htbp]
\centering\small
\begin{tabular}{|c||c|c|c||c|c|c|}
\hline
\multirow{2}{*}{DOF}& \multicolumn{3}{c||}{Unit Square}&\multicolumn{3}{|c|}{Hexagon}\\
\cline{2-7} &$\lambda_{1}$&$\lambda_{N_S}$&$\kappa$&$\lambda_{1}$&$\lambda_{N_S}$&$\kappa$
\\
\hline 512&1.22&9.34&7.64&0.29&4.88&16.63\\
\hline 1024&1.07&9.34&8.72&0.25&4.90&19.73\\
\hline 2048&0.95&9.34&9.87&0.21&4.91&23.38\\
\hline
\end{tabular}
\vskip 0.2cm \caption{Extremal eigenvalues of $T_{h,2}S_h$: convex domains.}
\label{tab:expV}
\end{table}

\begin{table}[h!!]
\centering\small
\begin{tabular}{|c||c|c|c|c||c|c|c|c|c|}
\hline
\multirow{2}{*}{DOF}& \multicolumn{4}{c||}{``L-"shape Domain($m_0=1$)}&\multicolumn{5}{|c|}{Trident($m_0=2$)}\\
\cline{2-10} &$\lambda_1$&$\lambda_{N_S-1}$&$\lambda_{N_S}$&$\kappa^{\rm eff}_{m_0}$&$\lambda_1$&$\lambda_{N_S-2}$&$\lambda_{N_S-1}$&$\lambda_{N_S}$&$\kappa^{\rm eff}_{m_0}$
\\
\hline 512&0.33&3.25&65.57&9.75&0.19&3.94&22.07&24.66&20.33\\
\hline 1024&0.31&3.37&102.92&10.98&0.16&4.01&29.01&32.52&25.24\\
\hline 2048&0.28&3.49&162.00&12.35&0.13&4.05&38.07&42.79&30.62\\
\hline
\end{tabular}
\vskip 0.2cm \caption{Extremal eigenvalues of $T_{h,2}S_h$: nonconvex domains.}\label{tab:expVII}
\end{table}

\begin{table}[h!!]
\centering\small
\begin{tabular}{|c||c|c|c||c|c|c|}
\hline
\multirow{2}{*}{DOF}& \multicolumn{3}{c||}{Unit Square}&\multicolumn{3}{|c|}{Hexagon}\\
\cline{2-7} &$\lambda_{1}$&$\lambda_{N_S}$&$\kappa$&$\lambda_{1}$&$\lambda_{N_S}$&$\kappa$
\\
\hline 352&0.23&9.33&40.64&0.027&4.88&180\\
\hline 704&0.21&9.33&45.46&0.023&4.90&216\\
\hline 1408&0.18&9.33&50.48&0.020&4.91&252\\
\hline
\end{tabular}
\vskip 0.2cm \caption{Extremal eigenvalues of $T_{h,3}S_h$: convex domains.}\label{tab:expVI}
\end{table}

\begin{table}[h!!]
\centering\small
\begin{tabular}{|c||c|c|c|c||c|c|c|c|c|}
\hline
\multirow{2}{*}{DOF}& \multicolumn{4}{c||}{``L-"shape Domain($m_0=1$)}&\multicolumn{5}{|c|}{Trident($m_0=2$)}\\
\cline{2-10} & $\lambda_1$ &  $\lambda_{N_S-1}$ & $\lambda_{N_S}$ & $\kappa^{\rm eff}_{m_0}$ & $\lambda_1$ & $\lambda_{N_S-2}$ & $\lambda_{N_S-1}$ & $\lambda_{N_S}$ & $\kappa^{\rm eff}_{m_0}$
\\
\hline 448&0.067&3.25&35.19&48.53&0.035&3.94&10.20&11.00&114\\
\hline 896&0.059&3.37&54.08&57.04&0.030&4.01&12.71&13.83&136\\
\hline 1792&0.053&3.49&83.84&66.14&0.025&4.04&15.93&17.45&159\\
\hline
\end{tabular}
\vskip 0.2cm \caption{Extremal eigenvalues of and $T_{h,3}S_h$: nonconvex domains.}\label{tab:expVIII}
\end{table}

\section{Concluding remarks}\label{sec:con}
This paper is aimed at developing efficient iterative methods for solving the algebraic systems arising from direct finite element discretization of the boundary value problems of fourth-order equations on an unstructured grid. The following objectives have been accomplished:
\begin{itemize}
\item[-] A class of nearly optimal iterative methods are developed for discrete fourth-order problem with both the first and second kinds boundary value conditions in two dimensions.
\item[-] A complete and rigorous analysis is provided for all the algorithms proposed in the paper.
\item[-] Numerical experiments are carried out to confirm all the theoretical results in the paper.
\end{itemize}
The algorithms and theories in this paper are valid for general unstructured grids in general polygonal domains which can be both convex and nonconvex. The iterative algorithms developed are the first and the only known methods in the literature for fourth-order problems that are provably (nearly) optimal.

To accomplish the objectives, the ÒFast Auxiliary Space PreconditioningÓ (FASP) method (\cite{Xu1996,Xu2010}) is used as the technical framework for designing the preconditioners, and the solution of a fourth-order problem is reduced to several second-order problems on the discrete level together with local relaxation methods. A number of intricate Sobolev spaces (such as the normal derivative trace space of $H^2(\Omega)\cap H^1_0(\Omega)$) defined on the boundary of a polygonal domains are carefully studied and thereafter used in the analysis of preconditioners. A special observation can be made that a straightforward mixed finite element discretization is used as the major component in the proposed preconditioners. Indeed, the aforementioned mixed method is either non-optimal or non-convergent method as a discretization method (\cite{Scholz1978,Grisvard1992,NazarovPlamenevsky1994,KozlovMazyaRossmann2001,ZhangZhang2008}) for the original fourth-order problem, and it is interesting to notice that the mixed method provides a nearly optimal preconditioner when it is used in conjuncture with additional local smoothers and preconditioned conjugate gradient methods.

Finally, the algorithms and theories in the paper need to be extended to the following cases in the future works:
\begin{itemize}
\item[-] General shape regular unstructured grids that are not assumed to be quasi-uniform.
\item[-] Three dimensional case.
\item[-] More complicated fourth-order equations such as Cahn--Hilliard equations.
\end{itemize}

\subsection*{Bibliographic comments}

The analysis in this paper consists of numerous technical results. Here we give a brief description of how some of these results are related to existing results in the literature.

The right inequality of \eqref{controltwoside} for first biharmonic problem was first studied in Babu\v{s}ka, Osborn, and Pitk\"aranta \cite{BabuskaOsbornPitkaranta1980} for convex domain in their analysis for some mixed methods, and Hanisch \cite{Hanisch2006} for nonconvex domain also for analysis of mixed methods. In the present paper, we establish \eqref{controltwoside} for both first and second biharmonic problems on both convex and nonconvex domains. 

Theorem \ref{lem:trace} can also be found in different form in Peisker \cite{Peisker1988}, who gave a proof assuming the domain is convex. In the present paper, we establish Theorem \ref{lem:trace} and prove it by means of an auxiliary Stokes problem on both convex and nonconvex domains. Lemma \ref{lem:gips} presents the isomorphisms between the normal derivative trace space and the Laplacian trace space of biharmonic functions on polygonal domains. When $\Omega$ is smooth, similar result was given by Glowinski and Pironneau \cite{GlowinskiPironneau1979}. 

Peisker \cite{Peisker1988} studied Lemma \ref{lem:bihesta} and Braess and Peisker \cite{BraessPeisker1986} proved Lemma \ref{lem:cprdeltas} both in a special case that $\Omega$ is convex and by analyzing the property of the discrete extension operator $E_h$ associated with the harmonic extension operator $E$. In the present paper, we study the property of $E_h$ associated with the generalized harmonic extension operator $E_c$, which does not coincide with $E$ in $H^{-1/2}_c(\Gamma)\setminus H^{1/2}(\Gamma)$ or in nonconvex domains, and establish Lemmas \ref{lem:bihesta} and \ref{lem:cprdeltas} for both convex and nonconvex domains. A disguised form of Lemma \ref{lem:deccrmix} can also be found in Glowinski and Pironneau \cite{GlowinskiPironneau1979}, but our formulation and proof are quite different.

Our preconditioner \eqref{FASPreconditioner} for the first biharmonic problem was motivated by an algorithm proposed in Peisker and Braess \cite{PeiskerBraess1987}, where an algebraic preconditioner for Morley element problem was presented for the special case that the domain is convex and the triangulation $\mathcal{T}_h$ consists of triangles that are similar to each other. Peisker \cite{Peisker1988} also noticed the role of $S_h$ and presented a preconditioner for $S_h$ in matrix form when $\Omega$ is convex. Her preconditioner can be realized by fast Fourier transform (FFT) on a graded bisection mesh of $\partial\Omega$. The technique of analyzing an interface operator by Fourier analysis was also used in, e.g., \cite{Chan.T;Sun.J;E.W1991,Dear.E;Glowinski.R;Pironneau.O1991}. 
 In the present paper, we design preconditioners in a unified framework for both first and second biharmonic problems discretized by various finite element methods on general triangulation for both convex and nonconvex domains. And, we establish preconditioners for $S_h$ in a more analytic approach, for general polygonal domain $\Omega$ that are triangulated by general meshes.

\end{document}